\documentclass[draft,preprint,12pt,numbers,sort&compress]{elsarticle}
\usepackage{mathrsfs}
\usepackage{amssymb}
\usepackage{amsthm}
\usepackage{mathrsfs}
\usepackage[centertags]{amsmath}
\usepackage{amsfonts}

\usepackage{color}

\input amssym.def
\input amssym.tex

\date{}

\newtheorem{Theorem}{Theorem}[section]

\newtheorem{Lemma}{Lemma}[section]

\newcommand\R{\mbox{\bf R}}
\newcommand\N{\mbox{\bf N}}
\newcommand\SR{\mbox{\scriptsize\bf R}}

\newcommand{\definition}{{\lower .5ex
  \hbox{$\>\>\stackrel{\triangle}{=}\>\>$} }}
\newcommand\supp{\mathop{\rm supp}}

\begin{document}

\baselineskip=22pt
\thispagestyle{empty}

\title{\huge The Cauchy problem for  a higher order shallow water  type equation on the circle\footnote{Corresponding author: W. YAN,\quad Email:\quad yanwei19821115@sina.cn}}

\author{Wei Yan\\
{\small Department of Mathematics and information science, Henan Normal University,}\\
{\small  Xinxiang, Henan\quad 453007, P.R.China} \\
{\small Email:yanwei19821115@sina.cn}\\
[2mm]
Yongsheng  Li\\
{\small Department of Mathematics, South China University of Technology}\\
{\small Email: Guangzhou, Guangdong 510640, P. R. China}\\
[2mm]
and\\
[2mm]
Jianhua Huang\\
{\small College of Science, National University of Defense and Technology,}\\
{\small Changsha,P. R. China\quad  410073}
}

\date{}
\maketitle

\noindent{\bf Abstract.}  In this paper, we investigate the Cauchy problem
for a higher order  shallow water type equation
\begin{eqnarray*}
     u_{t}-u_{txx}+\partial_{x}^{2j+1}u-\partial_{x}^{2j+3}u+3uu_{x}-2u_{x}u_{xx}-uu_{xxx}=0,
\end{eqnarray*}
 where  $x\in \mathbf{T}=\R/2\pi$ and $j\in N^{+}.$  Firstly, we prove that the Cauchy problem for the
 shallow water type equation
 is locally well-posed in $H^{s}(\mathbf{T})$ with $s\geq -\frac{j-2}{2}$ for arbitrary initial data.
 By using the $I$-method,  we prove that the Cauchy problem for the shallow water type equation is
 globally well-posed in $H^{s}(\mathbf{T})$ with $\frac{2j+1-j^{2}}{2j+1}<s\leq 1.$ Our results improve the result of A. A. Himonas, G. Misiolek (Communications in partial Differential Equations, 23(1998), 123-139;Journal of  Differential Equations, 161(2000), 479-495.)

\vskip 4mm
\noindent {\bf Keywords}: Periodic higher-order shallow water type  equation; Cauchy problem;
Low regularity

\vskip 4mm
\noindent {\bf Short Title:} Cauchy problem for a  shallow water type equation

\vskip 4mm
\noindent {\bf AMS  Subject Classification}:  35G25

\newpage

{\large\bf 1. Introduction}
\bigskip

\setcounter{Theorem}{0} \setcounter{Lemma}{0}

\setcounter{section}{1}

In this paper, we consider the Cauchy problem for a higher order shallow water  type equation
\begin{eqnarray}
&& u_{t}-u_{txx}+\partial_{x}^{2j+1}u-\partial_{x}^{2j+3}u+3uu_{x}-2u_{x}u_{xx}-uu_{xxx}= 0,\label{1.01}\\
    &&u(x,0)=u_{0}(x),\quad x\in \mathbf{T}=\R/2\pi,\label{1.02}
\end{eqnarray}
which is considered as the higher modification of the Camassa-Holm equation.
Rewrite (\ref{1.01})  as follows:
\begin{eqnarray}
    u_{t}+\partial_{x}^{2j+1}u+\frac{1}{2}\partial_{x}(u^{2})+
    \partial_{x}(1-\partial_{x}^{2})^{-1}\left[u^{2}+\frac{1}{2}u_{x}^{2}\right]=0,  \label{1.03}
\end{eqnarray}
which was derived by Camassa and Holm as a nonlinear model for water wave motion
in shallow channels with  the aid of  an asymptotic expansion directly in the Hamiltonian
for Euler equations \cite{CH,FF}.
 Omitting  the  last term yields
 \begin{eqnarray}
    u_{t}+\partial_{x}^{2j+1}u+\frac{1}{2}\partial_{x}(u^{2})=0.  \label{1.04}
\end{eqnarray}
When $j=1,$
equation (\ref{1.01})  reduces to the Korteweg-de Vries (KdV) equation
\begin{eqnarray}
u_{t}+u_{xxx}+\frac{1}{2}\partial_{x}(u^{2})=0.\label{1.05}
\end{eqnarray}
Kenig et. al. \cite{KPV, KPV0}  proved that   $s=-3/4$ is the critical Sobolev index for the KdV
equation in real line  and proved that the Cauchy problem for the periodic KdV equation is
locally well-posed in $H^{s}(0,2\pi\lambda)$ with $s\geq-\frac{1}{2}$ and $\lambda\geq1.$    Bourgain \cite{Bou}
proved that the Cauchy problem for the periodic KdV  equation is ill-posed in $H^{s}(0,2\pi\lambda)$
with $s< -\frac{1}{2}$ and $\lambda\geq 1.$  Colliander et.al. \cite{CKSTT} proved that the Cauchy problem for
the periodic KdV  equation is globally well-posed in $H^{s}(0,2\pi\lambda)$ with $s\geq -\frac{1}{2}$ and $\lambda\geq1.$
Kappeler and Topalov \cite{KT2005,KT2006} proved the global well-posedness of the
KdV and the defocusing mKV equations in $H^{s}(0,2\pi\lambda)$ for respectively $s\geq-1$
and $s\geq0$ and $\lambda\geq1$ with a solution-map which is continuous from $H^{-1}(0,2\pi\lambda) $(
$L^{2}(0,2\pi\lambda)$) into $C(R;H^{-1}(0,2\pi\lambda)$) ($C(R;L^{2}(0,2\pi\lambda))$) with $\lambda\geq1.$
Molinet \cite{M,Molinet} proved that the Cauchy problem for the periodic KdV equation
is ill-posed in $H^{s}(0,2\pi\lambda)$ with $s<-1$ and $\lambda\geq1$ in the  sense   that  the solution-map associated with the  KdV
equation is discontinuous for the $H^{s}(T)$ topology for  $s<-1$.

Lots of people have investigated the Cauchy problem for (\ref{1.03}), for instance, see \cite{By,CH,FF,HM1998,HM2000,HM1,LJ,O,LYY,LY,WC,YLL,YL}.
Himonas and Misiolek \cite{HM1998} proved that the Cauchy problem for (\ref{1.01})  is  locally well-posed
 for small initial data in $H^{s}(\mathbf{T})$ with $s\geq \frac{2-j}{2}$ and globally well-posed
 in $H^{1}(\mathbf{T})$.
 Himonas and Misiolek \cite{HM2000} proved that the Cauchy problem for (\ref{1.01}) with $j=1$ is
 locally well-posed
 for arbitrary initial data in $H^{s}(\mathbf{T})$ with $s\geq \frac{2-j}{2}$ and globally well-posed
 in $H^{1}(\mathbf{T})$.
Gorsky \cite{G} proved that the Cauchy problem for (\ref{1.01}) with $j=1$ is locally well-posed in
$H^{1/2}(\mathbf{T})$ for small initial data.
Li and Yang \cite{LY} prove that the Cauchy problem for (\ref{1.01}) with $j=1$ is locally well-posed
in $H^{s}(\mathbf{T})$ for $\frac{1}{2}<s<1$ and globally well-posed in
in $H^{s}(\mathbf{T})$ for $\frac{2}{3}<s<1$ with the aid of $I$-method. Olson \cite{O} proved that
the Cauchy problem for (\ref{1.01})  is locally well-posed in $H^{s}(\R)$ with $s>s^{'}$, where
$\frac{1}{4}\leq s^{'}<\frac{1}{2}.$ Yan et.al \cite{LYY} prove that the Cauchy problem for (\ref{1.01})
is locally well-posed in $H^{s}(\R)$ with $s>-j+\frac{5}{4}$ and is globally well-posed in $H^{1}(\R)$.
Yan et. al \cite{YLL} prove that  the Cauchy problem for (\ref{1.01})  is locally well-posed in
$H^{s}(\R)$ with $s=-j+\frac{5}{4}$,  $j\geq 2,j\in N^{+}$ and  ill-posed in $\dot{H}^{s}(\R)$
with $s<-j+\frac{5}{4}.$

In this paper, by establishing
some bilinear estimates and the fixed point Theorem, we prove that
the Cauchy problem for (\ref{1.01}) is locally well-posed in $H^{s}(\mathbf{T})$ with $s\geq \frac{2-j}{2}$;
by using the $I$-method, we prove that the problem is globally well-posed in $H^{s}(\mathbf{T})$ with $\frac{2j+1-j^{2}}{2j+1}<s\leq 1$.

We give some notations before stating the main results.  $0<\epsilon<\frac{1}{10000(2j+1)}$ and $\epsilon^{'}=\frac{1}{100(2j+1)}.$ $C$ is a positive constant
 which may vary from line to line.  $A\sim B$ means that $|B|\leq |A|\leq 4|B|$.
 $A\gg B$ means that $|A|\geq 4|B|.$ $a\vee b={\rm max}\left\{a,b\right\}.$ $a\wedge b={\rm min}\left\{a,b\right\}.$Let $\eta(t)$
 the smooth function supported in $[-1,2]$ and equals to
 $1$ in $[0,1]$.   Let $\Psi \in C_{0}^{\infty}(\R)$ be an even function such
  that $\Psi \geq 0,$ $\supp \Psi \subset [-\frac{3}{2},\frac{3}{2}]$,
 $\Psi= 1$ on $[-\frac{5}{4},\frac{5}{4}]$ and
 $v_{k}=\Psi(2^{-k}\xi)-\Psi(2^{-k+1}\xi).$

 For $k=k_{1}+k_{2},$  we  define
  \begin{eqnarray*}
|k_{min}|=\left\{|k|,|k_{1}|,|k_{2}|\right\},\quad |k_{max}|=\left\{|k|,|k_{1}|,|k_{2}|\right\}.
\end{eqnarray*}
 Throughout this paper,
 $\dot{Z}:=Z- \{ 0\}$ and $\dot{Z}^{+}:=Z^{+}- \{ 0\}$.    Denote
 $dk$ by the normalized counting measure on $\dot{Z}$:
 \begin{eqnarray*}
 \int a(k)dk=\sum_{k\in  \dot{Z}}a(k).
 \end{eqnarray*}
 Denote $\mathscr{F}_{x}f$ by the Fourier transformation of a function $f$
 defined on
  $[0,2\pi]$ with the respect to the space variable
 \begin{eqnarray*}
 \mathscr{F}_{x}f(k)=\frac{1}{2\pi}\int_{0}^{2\pi}e^{- i kx}f(x)dx.
 \end{eqnarray*}
 and we have the Fourier inverse transformation formula
 \begin{eqnarray*}
 f(x)=\int e^{ i kx} \mathscr{F}_{x}f(k)dk=\sum_{k \in
 \dot{Z}}e^{ i kx}\mathscr{F}_{x}f(k).
 \end{eqnarray*}
 Denote $\mathscr{F}_{t}f$ by the Fourier transformation of a function $f$
  with the respect to the time variable
 \begin{eqnarray*}
 \mathscr{F}_{t}f(\tau)=\frac{1}{2\pi}\int_{\SR}e^{- i t\tau}f(t)dt.
 \end{eqnarray*}
 and we have the Fourier inverse transformation formula
 \begin{eqnarray*}
 f(t)=\int e^{ i t\tau} \mathscr{F}_{t}f(\tau)d\tau.
 \end{eqnarray*}
 We define
\begin{eqnarray*}
S(t)\phi(x)=\int e^{i kx }e^{it k^{2j+1}}\mathscr{F}_{x}\phi(k)dk.
\end{eqnarray*}
We define the space-time  Fourier transform
$\mathscr{F}f(k,\tau)$ for $k\in \dot{Z}$ and $\tau\in \R$ by
\begin{eqnarray*}
\mathscr{F}f(k,\tau)=\frac{1}{2\pi}\int\int_{0}^{2\pi}e^{-i kx}e^{-i \tau t}f(x,t)dxdt
\end{eqnarray*}
and this transformation is inverted by
\begin{eqnarray*}
v(x,t)=\int\int e^{i kx}e^{i \tau t}\mathscr{F}f(k,\tau)dkd\tau.
\end{eqnarray*}
We define
\begin{eqnarray*}
\mathscr{F}_{x}\left[J_{x}^{s}\phi\right](k)=\langle k\rangle^{s}\mathscr{F}_{x}\phi (k),\mathscr{F}_{t}\left[J_{t}^{s}\phi\right](\tau)
=\langle \tau\rangle^{s}\mathscr{F}_{x}\phi (\tau).
\end{eqnarray*}
Thus, by using the above definitions, we have that
\begin{eqnarray*}
&&\|f\|_{L^{2}([0,2\pi ])}=\|\mathscr{F}_{x}f\|_{L^{2}(dk)},\\
&&\int_{0}^{2\pi }f(x)\overline{g(x)}dx=\int \mathscr{F}_{x}f(k)\overline{\mathscr{F}_{x}
f(k)}dk,\\
&&\mathscr{F}_{x}(fg)=\mathscr{F}_{x}f*\mathscr{F}_{x}g=\int \mathscr{F}_{x}f(k-k_{1})
\mathscr{F}_{x}g(k_{1})dk_{1}.
\end{eqnarray*}
Let
\begin{eqnarray*}
&&P(k)=k^{2j+1},\sigma=\tau+P(k),\quad \sigma_{l}=\tau_{l}+P(k_{l}),\quad l=1,2.
\end{eqnarray*}
For $s<1,$ we define
\begin{eqnarray*}
\mathscr{F}_{x}Iu(k)=m(k)\mathscr{F}_{x}u(k),
\end{eqnarray*}
where
$m(k)=\left(\frac{|k|}{N}\right)^{1-s}$ if $|k|>2N$,
$m(k)=1$ if $|k|\leq N$.
We define the Sobolev space $H^{s}(0,2\pi)$ with the norm
\begin{eqnarray*}
\|f\|_{H^{s}(\mathbf{T})}=\|\mathscr{F}_{x}f(k)\langle k\rangle^{s}\|_{L^{2}(k)}
\end{eqnarray*}
and define the $X_{s,\>b}$ spaces for $2\pi$-periodic KdV via the norm
\begin{eqnarray*}
\|u\|_{X_{s,b}(\mathbf{T}\times \SR)}=\left\|\langle k\rangle^{s} \left\langle \tau
+P(k)\right\rangle^{b}\mathscr{F}u(k,\tau)\right\|_{L^{2}(k\tau)}.
\end{eqnarray*}
and define the $Y_{s}$ space defined via the norm
\begin{eqnarray*}
\|u\|_{Y_{s}}=\|u\|_{X_{s,\frac{1}{2}}}+\left\|\langle k\rangle^{s}\mathscr{F}
u(k,\tau)\right\|_{L^{2}(k)L^{1}(\tau)}
\end{eqnarray*}
and define the $Z_{s}$ space defined via the norm
\begin{eqnarray*}
\|u\|_{Z_{s}}=\|u\|_{X_{s,-\frac{1}{2}}}+\left\|\frac{\langle k\rangle^{s}
\mathscr{F}u(k,\tau)}{\left\langle \tau+P(k)\right\rangle^{1/2}}\right\|_{L^{2}(k)L^{1}(\tau)}.
\end{eqnarray*}
We define
\begin{eqnarray*}
&&\|u\|_{X_{s,b}^{\delta}}={\rm inf}\left\{\|v\|_{X_{s,\>b}}\qquad v|_{[0,\>\delta]}=u\right\},\nonumber\\
&&\|u\|_{Y_{s}^{\delta}}={\rm inf}\left\{\|v\|_{Y_{s}}\qquad v|_{[0,\>\delta]}=u\right\}.
\end{eqnarray*}

The main result of this paper are as follows.

\begin{Theorem}\label{Thm1}
Let $s\geq -\frac{j-2}{2}$  and $u_{0}$ be $2\pi $-periodic
function and  zero $x$-mean.
Then the Cauchy problems (\ref{1.01})(\ref{1.02})
are locally well-posed in $H^{s}(\mathbf{T})$.
\end{Theorem}
\begin{Theorem}\label{Thm2}
Let  $\frac{2j+1-j^{2}}{2j+1}<s\leq 1$   and $u_{0}$ be $2\pi $-periodic
function and   zero $x$-mean.
Then the Cauchy problem (\ref{1.01})(\ref{1.02})
is globally well-posed in $H^{s}(\mathbf{T})$. More precisely, for any $T>0,$ let $u_{0}$ be $2\pi $-periodic
function and   zero $x$-mean, then  the Cauchy problems (\ref{1.01})(\ref{1.02})
are globally well-posed on $[0,T]$ in $H^{s}(\mathbf{T})$ with $\frac{2j+1-j^{2}}{2j+1}<s\leq 1$. Moreover,
\begin{eqnarray}
\sup\limits_{t\in [0,T]}\|u(\cdot,t)\|_{H^{s}}\leq
 CT^{\frac{1-s}{j-f(j)(1-s)}}\|u_{0}\|_{H^{s}}^{\frac{j}{j-f(j)(1-s)}},
\end{eqnarray}
where
\begin{eqnarray*}
f(j)=\frac{(2j+1)}{j-3(2j+1)\epsilon}.
\end{eqnarray*}
\end{Theorem}

The rest of the paper is arranged as follows. In Section 2,  we give some
preliminaries. In Section 3, we establish the bilinear estimate. In Section 4, we give the proof of Theorem 1.1.
In Section 5, we give the proof of Theorem 1.2.

\bigskip
\bigskip
\bigskip

 \noindent{\large\bf 2. Preliminaries }

\setcounter{equation}{0}

\setcounter{Theorem}{0}

\setcounter{Lemma}{0}

\setcounter{section}{2}
In this section, we make some preliminaries which are crucial in establishing the Theorem 1.1.

\begin{Lemma}\label{Lemma2.1}
Let $u_{l}$ with $l=1,2$ be $L^{2}(\dot{Z}\times \R)$-real valued functions. Then for any $(l_{1},l_{2}) \in \N^{2}$
\begin{eqnarray}
      \left\|(\Psi_{l_{1}}u_{1})*(\Psi_{l_{2}}u_{2})\right\|_{L_{xt}^{2}}\leq C\left(2^{l_{1}}
      \wedge2^{l_{2}}\right)^{1/2}\left(2^{l_{1}}\vee2^{l_{2}}\right)^{\frac{1}{2(2j+1)}}
      \|\Psi_{l_{1}}u_{1}\|_{L^{2}}\|\Psi_{l_{2}}u_{2}\|_{L^{2}}.
        \label{2.01}
\end{eqnarray}
\end{Lemma}
{\bf Proof.} As the proof of \cite{Bou,ST}, we can assume that $\supp u_{l}
\subset \left\{(\tau,k)\in \R\times \dot{Z}^{+}\right\}$.
By using the Cauchy-Schwarz in $(\tau_{1},k_{1})$, we have that
\begin{eqnarray}
&&\left\|(\Psi_{l_{1}}u_{1})*(\Psi_{l_{2}}u_{2})\right\|_{L^{2}}^{2}\nonumber\\&&
=\int_{\SR_{\tau}}\sum_{k\in \dot{Z}}\left|\int_{\SR_{\tau_{1}}}\sum_{k_{1}\in\dot{Z}}(\Psi_{l_{1}}u_{1})(\tau_{1},k_{1})
(\Psi_{l_{2}}u_{2})(\tau-\tau_{1},k-k_{1})d\tau_{1}\right|^{2}d\tau\nonumber\\&&\leq
C\int_{\tau}\sum_{k \in \dot{Z}}\alpha(\tau,k)\int_{\SR_{\tau_{1}}}
\sum_{k_{1}\in \dot{Z}}\left|(\Psi_{l_{1}}u_{1})(\tau_{1},k_{1})
(\Psi_{l_{2}}u_{2})(\tau-\tau_{1},k-k_{1})\right|^{2}d\tau_{1}d\tau\nonumber\\&&
\leq C\sup_{\tau\in \SR,\>k \in \dot{Z}}\alpha(\tau,k)
\|\Psi_{l_{1}}u_{1}\|_{L^{2}}^{2}\|\Psi_{l_{2}}u_{2}\|_{L^{2}}^{2},\label{2.02}
\end{eqnarray}
where
\begin{eqnarray*}
&&\hspace{-0.9cm}\alpha(\tau,k)\leq C\#\Lambda_{1}(\tau,k),
\end{eqnarray*}
 here
 \begin{eqnarray*}
\Lambda_{1}(\tau,k)= \left\{(\tau_{1},k_{1})\in \R\times
\dot{Z}^{+}/k-k_{1}\in \dot{Z}^{+}, \left\langle\sigma_{1}\right\rangle\sim 2^{l_{1}},
\left\langle\sigma_{2}\right\rangle\sim 2^{l_{2}}\right\}
 \end{eqnarray*}
 For fixed $\tau,\xi\neq 0$,  We define $M^{'}=\tau+(-1)^{j}\frac{\xi^{2j+1}}{4^{j}}$  and
 let $E_1$ and $E_2$ be the projections of $\Lambda_{1}$ onto
the $k_1$-axis and $\tau_1$-axis, respectively.
It is easily checked  that
\begin{eqnarray}
&&\left(\tau+(-1)^{n}\frac{k^{2j+1}}{4^{j}}\right)-
\left(\tau_{1}+(-1)^{j}k_{1}^{2j+1}\right)-\left(\tau_{2}+(-1)^{j}k_{2}^{2j+1}\right)\nonumber\\&&
=(-1)^{j+1}\left[k_{1}^{2j+1}+k_{2}^{2j+1}-\frac{k^{2j+1}}{4^{j}}\right]=(-1)^{j+1}
k\left(k_{1}-\frac{k}{2}\right)^{2}F(k,k_{1}),\label{2.03}
\end{eqnarray}
where
\begin{eqnarray*}
&&F(k,k_{1})\nonumber\\&&\hspace{-1cm}=C_{2j+1}^{2}\left(\frac{1}{2}\right)^{2j-2}k^{2j-2}+ C_{2j+1}^{4}
\left(\frac{1}{2}\right)^{2j-2}k^{2j-4}\left(k_{1}-\frac{k}{2}\right)^{2}+
\cdot\cdot\cdot +C_{2j+1}^{2j}\left(k_{1}-\frac{k}{2}\right)^{2j-2}.
\end{eqnarray*}
From (\ref{2.03}), we have that there exist two constant $C_{1},C_{2}>0$ such that
\begin{eqnarray}
\frac{\left|C_{1}(2^{l_{1}}+ 2^{l_{2}})-M^{'}\right|}{|kF(k,k_{1})|}\leq
 \frac{3}{4}(k_{1}-k_{2})^{2}\leq\frac{\left| C_{2}(2^{l_{1}}+ 2^{l_{2}})+M^{'}\right|}{|k F(k,k_{1})|},\label{2.04}
\end{eqnarray}
When $k^{2j+1}>2^{l_{1}}\vee 2^{l_{2}},$ from (\ref{2.04}), we have that
\begin{eqnarray}
&&\#E_{2}\leq {\rm mes \> E_{2}}+1\leq 2\left[\frac{\left|
 C_{1}(2^{l_{1}}+ 2^{l_{2}})+M^{'}\right|}{|k F(k,k_{1})|}-
 \frac{\left|C_{2}(2^{l_{1}}+ 2^{l_{2}})-M^{'}\right|}{|k F(k,k_{1})|}\right]^{1/2}+1\nonumber\\&&
 \leq C\left(\frac{ (2^{l_{1}}\vee 2^{l_{2}})}{|k^{2j-1}|}\right)^{1/2}+1
 \leq C \left(2^{l_{1}}\vee 2^{l_{2}}\right)^{\frac{1}{2j+1}} \label{2.05}.
\end{eqnarray}
When $0\leq k^{2j+1}\leq 2^{l_{1}}\vee 2^{l_{2}}$, since $0\leq k_{1}\leq k,$ we have that
\begin{eqnarray}
&&\#E_{2}\leq  \#\left\{k_{1},\quad 0\leq k_{1}^{2j+1}\leq 2^{l_{1}}
\vee 2^{l_{2}}\right\}\leq C \left(2^{l_{1}}\vee 2^{l_{2}}\right)^{\frac{1}{2j+1}} \label{2.06}.
\end{eqnarray}
From (\ref{2.02}), it is easily checked that
\begin{eqnarray}
\#E_{1}\leq {\rm mes}\>E_{1}+1\leq C \left(2^{l_{1}}\wedge 2^{l_{2}}\right).\label{2.07}
\end{eqnarray}
Combining (\ref{2.02})  with  (\ref{2.05})-(\ref{2.07}),  we have that
\begin{eqnarray}
&&\left\|(\Psi_{l_{1}}u_{1})*(\Psi_{l_{2}}u_{2})\right\|_{L^{2}}\leq C\left(2^{l_{1}}
\wedge 2^{l_{2}}\right)^{1/2}\left(2^{l_{1}}\vee 2^{l_{2}}\right)^{\frac{1}{2(2j+1)}}\|\Psi_{l_{1}}u_{1}\|_{L^{2}}\|\Psi_{l_{2}}u_{2}\|_{L^{2}}.\label{2.08}
\end{eqnarray}

We have completed the proof of Lemma 2.1.

\begin{Lemma}\label{Lemma2.2}
Let $v(x,t)$ be a $2\pi $-periodic function. Then
\begin{eqnarray}
      \left\|v\right\|_{L_{xt}^{4}}\leq C\|v\|_{X_{0,\frac{(j+1)}{2(2j+1)}}(\mathbf{T} \times \SR)}.
        \label{2.09}
\end{eqnarray}
\end{Lemma}
{\bf Proof.} By using the triangle inequality, let $l_{1}=l+l_{2}$ with $l\in N,$ by using (\ref{2.01}),  we have that
\begin{eqnarray}
&&\|v\|_{L_{xt}^{4}}^{2}=\|v^{2}\|_{L^{2}}=\|\mathscr{F}v*\mathscr{F}v\|_{L^{2}}\leq
\sum\limits_{l_{1}\geq 0}\sum\limits_{l_{2}\geq 0}\left\|\Psi_{l_{1}}|\mathscr{F}v|\Psi_{l_{2}}|\mathscr{F}v|\right\|_{L^{2}}\nonumber\\
&&\leq C \sum\limits_{l_{1}\geq 0}\sum\limits_{l_{2}\geq 0}\left\|\Psi_{l_{1}}|\mathscr{F}v|*\Psi_{l_{2}}|\mathscr{F}v|\right\|_{L^{2}}
\nonumber\\&&\leq C\sum\limits_{l\geq 0}\sum\limits_{l_{2}\geq 0}2^{l_{2}/2}2^{(l_{2}+l)/2(2j+1)}\left\|\Psi_{l_{2}+l}
\mathscr{F}v\right\|_{L^{2}}\|\Psi_{l_{2}}\mathscr{F}v\|_{L^{2}}\nonumber\\&&
\leq C\sum\limits_{l\geq 0}\sum\limits_{l_{2}\geq 0}2^{\frac{j+1}{2(2j+1)}l_{2}}\|\Psi_{l_{2}}\mathscr{F}v\|_{L^{2}}2^{-\frac{j}{2(2j+1)}l}2^{\frac{(j+1)(l_{2}+l)}{2(2j+1)}}
\left\|\Psi_{l_{2}+l}\mathscr{F}v\right\|_{L^{2}}\nonumber\\&&
\leq C\sum\limits_{l\geq 0}2^{-\frac{j}{2(2j+1)}l}\left(\sum_{l_{2}\geq 0}2^{\frac{j+1}{2(2j+1)}l_{2}}\|\Psi_{l_{2}}\mathscr{F}v\|_{L^{2}}^{2}\right)^{1/2} \left(2^{\frac{(j+1)(l_{2}+l)}{2(2j+1)}}\left\|\Psi_{l_{2}+l}\mathscr{F}v\right\|_{L^{2}}^{2}\right)^{1/2}\nonumber\\
&&\leq C\|v\|_{X_{0,\frac{(j+1)}{2(2j+1)}}([0,\>2\pi ] \times \SR)}^{2}.\label{2.010}
\end{eqnarray}
From (\ref{2.010}),  we have (\ref{2.09}).

We have completed the proof of Lemma 2.2.

{\bf Remark:}  In  line -3 of page 493 in \cite{HM2000},  Himonas and Misiolek presented the conclusion of Lemma 2.2, however, the proof process  is not given.

\begin{Lemma}\label{Lemma2.3}
Let $v(x,t)$ be a $2\pi $-periodic function. Then
\begin{eqnarray}
      \left\|v\right\|_{X_{0,-\frac{(j+1)}{2(2j+1)}}(\mathbf{T} \times \SR)}
      \leq C\|v\|_{L_{xt}^{4/3}}=\left(\int_{0}^{2\pi} v^{4/3}(x,t)dxdt\right)^{3/4}.
        \label{2.011}
\end{eqnarray}
\end{Lemma}
{\bf Proof.}  Combining the Lemma 2.2 with the duality, we have Lemma 2.3.

\begin{Lemma}\label{Lemma2.4}Let
\begin{eqnarray*}
&&k=k_{1}+k_{2},\tau=\tau_{1}+\tau_{2},\nonumber\\
&&\sigma=\tau+(-1)^{j}k^{2j+1},\sigma_{l}=\tau_{l}+(-1)^{j}k_{l}^{2j+1},l=1,2.
\end{eqnarray*}
Then
\begin{eqnarray*}
3{\rm max}\left\{|\sigma|,|\sigma_{1}|,|\sigma_{2}|\right\}\leq|\sigma-\sigma_{1}-\sigma_{2}|
=\left|k^{2j+1}-k_{1}^{2j+1}-k_{2}^{2j+1}\right|\sim |k_{min}||k_{max}|^{2j}.
\end{eqnarray*}
\end{Lemma}

For the proof of Lemma 2.4, we refer the readers to  Lemma 2.5 in \cite{YLL}.

\begin{Lemma}\label{Lemma2.5}
For $k \in \dot{Z}$, $k_{j}\in \dot{Z}(j=1,2)$ and dyadic $M\geq 1$ and  $\epsilon^{'}=\frac{1}{100(2j+1)},$ we have that
\begin{eqnarray}
&&{\rm mes}\left\{\mu \in \R: \quad |\mu|\sim M, \mu=k^{2j+1}-k_{1}^{2j+1}-k_{2}^{2j+1}
+O(\langle |k_{\rm min}|k_{\rm max}|^{2j}\rangle^{\epsilon^{'}})\right\}\nonumber\\&&\leq CM^{\frac{100j+1}{50(2j+1)}}. \label{2.012}
\end{eqnarray}
\end{Lemma}
{\bf Proof.} Without loss  of generality, we can assume that $|k_{1}|\geq |k_{2}|.$
 When $|k|\geq |k_{1}|$ which yields that $|k_{1}|\leq |k|\leq 2|k_{1}|$, from
\begin{eqnarray}
 \mu=k^{2j+1}-k_{1}^{2j+1}-k_{2}^{2j+1}
+O(\langle |k_{\rm min}|k_{\rm max}|^{2j}\rangle^{\epsilon^{'}}), \label{2.013}
\end{eqnarray}
we have that $C_{1}|k|^{2j}\leq|\mu|\leq C_{2}|k|^{2j+1}$ since $k_{1},k_{2}\in \dot{Z}.$
Thus, we have that $|\mu|\sim M\sim |k|^{p}$, $p\in [2j,2j+1].$
Thus, $|k_{1}^{2j-1}k_{2}|\sim M^{1-\frac{1}{p}}$, $p\in [2j,2j+1].$ Consequently, we have that
\begin{eqnarray}
&&{\rm mes}\left\{\mu \in \R: \quad |\mu|\sim M, \mu=k^{2j+1}-k_{1}^{2j+1}-k_{2}^{2j+1}
+O(\langle|k_{\rm min}||k_{\rm max}|^{2j}\rangle^{\epsilon^{'}})\right\}\nonumber\\&&\leq CM^{\frac{2j}{2j+1}}
M^{\epsilon^{'}}\leq CM^{\frac{200j+1}{100(2j+1)}}.\label{2.014}
\end{eqnarray}
When $|k_{1}|\geq|k|,$ from (\ref{2.013}), we have that $C_{1}|k_{1}|^{2j}\leq|\mu|\leq C_{2}|k_{1}|^{2j+1}$ since $k_{1},k_{2}\in \dot{Z}.$
Thus, we have that $|\mu|\sim M\sim |k_{1}|^{p}$, $p\in [2j,2j+1].$
Thus, $|k_{1}^{2j-1}k|\sim M^{1-\frac{1}{p}}$, $p\in [2j,2j+1].$ Consequently, we have that
\begin{eqnarray}
&&{\rm mes}\left\{\mu \in \R: \quad |\mu|\sim M, \mu=k^{2j+1}-k_{1}^{2j+1}-k_{2}^{2j+1}
+O(\langle|k_{\rm min}||k_{\rm max}|^{2j}\rangle^{\epsilon^{'}})\right\}\nonumber\\&&\leq CM^{\frac{2j}{2j+1}}
M^{\epsilon^{'}}\leq CM^{\frac{200j+1}{100(2j+1)}}.\label{2.015}
\end{eqnarray}

We have completed the proof of Lemma 2.5.

\begin{Lemma}\label{Lemma2.6}Let $\phi$ be $2\pi$-periodic function. Then
\begin{eqnarray}
\left\|\eta(t)S(t)\phi\right\|_{Y_{s}^{\delta}}\leq C\|\phi\|_{H^{s}}.\label{2.016}
\end{eqnarray}
\end{Lemma}
{\bf Proof.} To obtain (\ref{2.016}),  it suffices to prove that
\begin{eqnarray}
\left\|\eta(t)\eta\left(\frac{t}{\delta}\right) S(t)\phi\right\|_{Y_{s}}\leq C\|\phi\|_{H^{s}}\label{2.017}.
\end{eqnarray}
From Lemma 7.1 of \cite{CKSTT}, we have that
\begin{eqnarray}
\left\|\eta(t)\eta\left(\frac{t}{\delta}\right) S(t)\phi\right\|_{Y_{s}}\leq C\|\eta\left(\frac{t}{\delta}\right)\phi\|_{H^{s}}\leq C\|\phi\|_{H^{s}}.\label{2.018}
\end{eqnarray}

We have completed the Lemma 2.6.

\begin{Lemma}\label{Lemma2.7}Let $F$ be $2\pi$-periodic function. Then
\begin{eqnarray}
\left\|\eta(t) \int_{0}^{t}S(t-\tau)F(\tau)d\tau\right\|_{Y_{s}^{\delta}}\leq
C\|\eta\left(\frac{t}{\delta}\right) F\|_{Z_{s}}.\label{2.019}
\end{eqnarray}
\end{Lemma}
{\bf Proof.} To obtain (\ref{2.019}),  it suffices to prove that
\begin{eqnarray}
\left\|\eta(t) \eta\left(\frac{t}{\delta}\right)\int_{0}^{t}S(t-\tau)F(\tau)d\tau\right\|_{Y_{s}^{\delta}}
\leq C\left\|\eta\left(\frac{t}{\delta}\right) F\right\|_{Z_{s}}\label{2.020}
\end{eqnarray}
which follows from Lemma 7.2 of \cite{CKSTT}.

We have completed the proof of Lemma 2.7.

\begin{Lemma}\label{Lemma2.8}
Let
\begin{eqnarray*}
\Omega(k)=\left\{\mu \in \R: \quad |\mu|\sim M, \mu=k^{2j+1}-k_{1}^{2j+1}-k_{2}^{2j+1}
+O(\langle |k_{\rm min}|k_{\rm max}|^{2j}\rangle^{\epsilon})\right\}
\end{eqnarray*}
 Then
\begin{eqnarray}
\int \langle\mu\rangle^{-1}\chi_{\Omega(k)}(\mu)d\mu\leq C.\label{2.021}
\end{eqnarray}
\end{Lemma}
{\bf Proof.} Combining Lemma 2.6 with the proof of
 page 737 in \cite{CKSTT},  we have Lemma 2.10.

\begin{Lemma}\label{Lemma2.9}Let $s\in \R$ and $\delta \in (0,1),$
then for
$-\frac{1}{2}<b<b^{'}\leq 0$ or $0\leq b<b^{'}<\frac{1}{2}$,
we have that
\begin{eqnarray}
&&\left\|\eta\left(\frac{t}{\delta}\right)u\right\|_{X_{0,\>b}}\leq
C\delta^{b-b^{'}}\|u\|_{X_{0,\>b^{'}}},\label{2.022}
\end{eqnarray}
\end{Lemma}

For the proof of Lemma 2.9, we refer the readers to Lemma  1.10
of \cite{G}.

\begin{Lemma}\label{Lemma2.10}
For $u\in X_{\sigma, b}^{\delta}$ there exists $\tilde{u}$ with $u|_{[0,\delta]}=\tilde{u}$, such that for $s\leq\sigma$, we have that
\begin{eqnarray*}
\|u\|_{X_{s,b}^{\delta}}=\|\tilde{u}\|_{X_{s,b}}.
\end{eqnarray*}
\end{Lemma}

For the proof of Lemma 2.10, we refer  the readers to Lemma 1.6 of \cite{G}.

\begin{Lemma}\label{Lemma2.11} Let $s\in \R$ and $0<\epsilon<\frac{1}{10000(2j+1)}$ and
\begin{eqnarray}
F(k,\tau)=\langle k\rangle^{s}\langle\sigma\rangle^{1/2}
\mathscr{F}\left(\eta\left(\frac{t}{\delta}\right)\tilde{u}\right)(k,\tau),\label{2.023}
\end{eqnarray}
where
$F\in L^{2}.$
Then
\begin{eqnarray}
\left\|\mathscr{F}^{-1}\left(\frac{F}{\langle\sigma\rangle^{1/2}}\right)\right\|_{L^{4}}\leq C\delta^{\frac{j}{2(2j+1)}-\epsilon}\|F\|_{L^{2}}.\label{2.024}
\end{eqnarray}
\end{Lemma}
{\bf Proof.} From (\ref{2.023}) and Lemmas 2.2, 2.9, we have that
\begin{eqnarray}
&&\left\|\mathscr{F}^{-1}\left(\frac{F}{\langle\sigma\rangle^{1/2}}\right)\right\|_{L^{4}}
=\left\|\eta\left(\frac{t}{\delta}\right)J_{x}^{s}\tilde{u}\right\|_{L^{4}}\nonumber\\&&
\leq C\left\|\eta\left(\frac{t}{\delta}\right)J_{x}^{s}\tilde{u}\right\|_{X_{0,\frac{j+1}{2(2j+1)}}}\nonumber\\
&&\leq C\delta^{\frac{j}{2(2j+1)}-\epsilon}\left\|\eta\left(\frac{t}{\delta}\right)J_{x}^{s}\tilde{u}\right\|_{X_{0,\frac{1}{2}-\epsilon}}\nonumber\\
&&\leq C\delta^{\frac{j}{2(2j+1)}-\epsilon}\left\|\eta\left(\frac{t}{\delta}\right)\tilde{u}\right\|_{X_{s,\frac{1}{2}}}\nonumber\\
&&=C\delta^{\frac{j}{2(2j+1)}-\epsilon}\|F\|_{L^{2}}.\label{2.025}
\end{eqnarray}

We have completed the proof of Lemma 2.11.

{\bf Remark:}  Lemma 2.11 improves the result of Lemma 3.2 in \cite{HM2000} with $\mu=2j+1$.

\begin{Lemma}\label{Lemma2.12} Let
\begin{eqnarray*}
\sigma=\tau +(-1)^{j}k^{2j+1},\sigma_{l}=\tau_{l}+(-1)^{j}k_{l}^{2j+1},l=1,2.
\end{eqnarray*}
 and  $s\in \R$ and $0<\epsilon<\frac{1}{10000(2j+1)}$ and
\begin{eqnarray}
G_{l}(k_{l},\tau_{l})=\langle k_{l}\rangle^{s}\langle\sigma_{l}\rangle^{1/2}
\mathscr{F}\left(\eta\left(\frac{t}{\delta}\right)\tilde{u}_{l}\right)(k_{l},\tau_{l}),l=1,2,\label{2.026}
\end{eqnarray}
where
$G_{l}\in L^{2},l=1,2.$
Then
\begin{eqnarray}
\left\|\langle\sigma\rangle^{-\frac{1}{2}+\epsilon}\int_{\!\!\!\mbox{\scriptsize $
\begin{array}{l}
k=k_{1}+k_{2}\\
\tau=\tau_{1}+\tau_{2}
\end{array}
$}}\frac{\prod_{l=1}^{2}G_{l}(k_{l},\tau_{l})}{\langle\sigma_{2}\rangle^{1/2}}dk_{1}
d\tau_{1} \right\|_{L^{2}}\leq C\delta^{\frac{j}{2j+1}-2\epsilon}\prod_{l=1}^{2}\|G_{l}\|_{L^{2}}.\label{2.027}
\end{eqnarray}
\end{Lemma}
{\bf Proof.} By using Lemmas 2.3, 2.4,2.11, we have that
\begin{eqnarray*}
&&\left\|\langle\sigma\rangle^{-\frac{1}{2}+\epsilon}\int_{\!\!\!\mbox{\scriptsize $
\begin{array}{l}
k=k_{1}+k_{2}\\
\tau=\tau_{1}+\tau_{2}
\end{array}
$}}\frac{\prod_{l=1}^{2}G_{l}(k_{l},\tau_{l})}{\langle\sigma_{2}\rangle^{1/2}}dk_{1}
d\tau_{1} \right\|_{L^{2}}\nonumber\\
&&=\left\|\left(\eta\left(\frac{t}{\delta}\right)J_{x}^{s}\tilde{u}_{2}\right)\mathscr{F}^{-1}(G_{1})\right\|_{X_{0,-\frac{1}{2}+\epsilon}}\nonumber\\
&&\leq C\delta^{\frac{j}{2(2j+1)}-\epsilon}\left\|\left(\eta\left(\frac{t}{\delta}\right)J_{x}^{s}\tilde{u}_{2}\right)\mathscr{F}^{-1}(G_{1})\right\|_{X_{0,-\frac{j+1}{2(2j+1)}}}\nonumber\\
&&\leq C\delta^{\frac{j}{2(2j+1)}-\epsilon}\left\|\left(\eta\left(\frac{t}{\delta}\right)J_{x}^{s}\tilde{u}_{2}\right)\mathscr{F}^{-1}(G_{1})\right\|_{L_{xt}^{4/3}}\nonumber\\
&&\leq C\delta^{\frac{j}{2(2j+1)}-\epsilon}\left\|\mathscr{F}^{-1}(G_{1})\right\|_{L^{2}}\left\|J_{x}^{s}\eta\left(\frac{t}{\delta}\right)\tilde{u}_{2}\right\|_{L_{xt}^{4}}\nonumber\\
&&\leq C\delta^{\frac{j}{2(2j+1)}-\epsilon}\left\|\mathscr{F}^{-1}(G_{1})\right\|_{L^{2}}\left\|J_{x}^{s}\eta\left(\frac{t}{\delta}\right)\tilde{u}_{2}\right\|_{X_{0,\frac{j+1}{2(2j+1)}}}\nonumber\\
&&\leq C\delta^{\frac{j}{2j+1}-2\epsilon}\left\|\mathscr{F}^{-1}(G_{1})\right\|_{L^{2}}\left\|J_{x}^{s}\eta\left(\frac{t}{\delta}\right)\tilde{u}_{2}\right\|_{X_{0,\frac{1}{2}-\epsilon}}\nonumber\\
&&\leq C\delta^{\frac{j}{2j+1}-2\epsilon}\left\|\mathscr{F}^{-1}(G_{1})\right\|_{L^{2}}\left\|J_{x}^{s}\eta\left(\frac{t}{\delta}\right)\tilde{u}_{2}\right\|_{X_{0,\frac{1}{2}}}\nonumber\\
&&\leq C\delta^{\frac{j}{2j+1}-2\epsilon}\prod_{l=1}^{2}\|G_{l}\|_{L^{2}}.
\end{eqnarray*}

We have completed the proof of Lemma 2.12.

\bigskip
\bigskip

\noindent{\large\bf 3. Bilinear estimates }

\setcounter{equation}{0}

 \setcounter{Theorem}{0}

\setcounter{Lemma}{0}

 \setcounter{section}{3}
In this section, we establish some important bilinear estimates which are the core of this paper

\begin{Lemma}\label{Lemma3.1}
Let $u_{l}(x,t)$ with $l=1,2$ which are zero $x$-mean for all $t$ be
$2\pi$- periodic functions  of $x$
and $s\geq \frac{2-j}{2}$. For $\epsilon<\frac{1}{10000(2j+1)},$
then we have that
\begin{eqnarray}
      \left\|\partial_{x}(1-\partial_{x}^{2})^{-1}\left[
      \prod_{l=1}^{2}\left[\partial_{x}\eta\left(\frac{t}{\delta}\right)u_{l}\right]\right]\right\|_{X_{s,-\frac{1}{2}}^{\delta}}\leq
      C\delta
      ^{\frac{j}{2j+1}-2\epsilon}
      \prod_{l=1}^{2}\|u_{l}\|_{X_{s,\frac{1}{2}}^{\delta}}.
        \label{3.01}
\end{eqnarray}
\end{Lemma}
{\bf Proof.}Let $\tilde{u}$ and $\tilde{u}_{1},\tilde{u}_{2}$ be
the extension of $u,u_{1},u_{2}$,
 respectively, according to Lemma 2.10, we have that
\begin{eqnarray*}
\|u\|_{X_{s,\frac{1}{2}}^{\delta}}=\|\tilde{u}\|_{X_{s,\frac{1}{2}}},
\|u_{l}\|_{X_{s,\frac{1}{2}}^{\delta}}
=\|\tilde{u}_{l}\|_{X_{s,\frac{1}{2}}},\quad l=1,2.
\end{eqnarray*}
By duality and  the Plancherel identity,
for $u\in X_{-s,\frac{1}{2}}^{\delta},$ to obtain (\ref{3.01}),
it suffices to prove that
\begin{eqnarray}
&&\int_{\!\!\!\mbox{\scriptsize $
\begin{array}{l}
k=k_{1}+k_{2}\\
\tau=\tau_{1}+\tau_{2}
\end{array}
$}}\left|\frac{kk_{1}k_{2}}{1+k^{2}}\mathscr{F}\left(\eta\left(\frac{t}{\delta}\right)\tilde{u}\right)
(k,\tau)\prod_{l=1}^{2}
\mathscr{F}\left(\eta\left(\frac{t}{\delta}\right)\tilde{u}_{l}\right)(k_{l},\tau_{l})
\right|dk_{1}
d\tau_{1}dkd\tau \nonumber\\&&\leq C\delta
      ^{\frac{j}{2j+1}-2\epsilon}
\|u\|_{X_{-s,\frac{1}{2}}^{\delta}}
\prod_{l=1}^{2}\|u_{l}\|_{X_{s,\frac{1}{2}}^{\delta}}=C\delta
      ^{\frac{j}{2j+1}-2\epsilon}\|\tilde{u}\|_{X_{-s,\frac{1}{2}}}\prod_{l=1}^{2}
      \|\tilde{u}_{l}\|_{X_{s,\frac{1}{2}}}.\label{3.02}
\end{eqnarray}
Without loss of generality, we can assume that
 $\mathscr{F}\left(\eta\left(\frac{t}{\delta}\right)\tilde{u}_{l}\right)
 (k_{l},\tau_{l})\geq 0(l=1,2)$
 and $\mathscr{F}\left(\eta\left(\frac{t}{\delta}\right)\tilde{u}\right)
 (k,\tau)\geq 0$.
Let
\begin{eqnarray*}
&&F(k,\tau)=\langle k\rangle^{-s}\langle\sigma\rangle^{1/2}\mathscr{F}
\left(\eta\left(\frac{t}{\delta}\right)\tilde{u}\right)(k,\tau),\\
&&F_{l}(k_{l},\tau_{l})=\langle k_{l}\rangle^{s}\langle\sigma_{l}\rangle^{1/2}
\mathscr{F}\left(\eta\left(\frac{t}{\delta}\right)\tilde{u}_{l}\right)(k_{l},\tau_{l})
,\quad l=1,2,\\
&&K_{1}(k_{1},\tau_{1},k,\tau)=\frac{|kk_{1}k_{2}|\langle k\rangle^{s}}{(1+k^{2})
\langle\sigma\rangle^{1/2}
\prod_{l=1}^{2}\langle k_{l}\rangle^{s}\langle\sigma_{l}\rangle^{1/2}}.\\
\end{eqnarray*}
To obtain (\ref{3.02}),  it suffices to prove that
\begin{eqnarray}
&&\int_{\!\!\!\mbox{\scriptsize $
\begin{array}{l}
k=k_{1}+k_{2}\\
\tau=\tau_{1}+\tau_{2}
\end{array}
$}}K_{1}(k_{1},\tau_{1},k,\tau)F(k,\tau)\prod_{l=1}^{2}F_{l}(k_{l},\tau_{l})
dk_{1}d\tau_{1} dkd\tau\nonumber
\\&&\leq C\delta
      ^{\frac{j}{2j+1}-2\epsilon}\|F\|_{L^{2}}\prod_{l=1}^{2}\|F_{l}\|_{L^{2}}\label{3.03}.
\end{eqnarray}
From the mean zero condition, we can assume that $k\neq0,k_{l}\neq0 (l=1,2).$

Since ${\rm min}\left\{|k|,|k_{1}|,|k_{2}|\right\}\geq 1,$ from Lemma 2.4,
we have that one of the following three cases must occur:
\begin{eqnarray}
&&(a):\quad |\sigma|={\rm max}\left\{|\sigma|,|\sigma_{1}|,|\sigma_{2}|\right\}
\geq C|k_{min}||k_{max}|^{2j},\nonumber\\
&&(b):\quad |\sigma_{1}|={\rm max}\left\{|\sigma|,|\sigma_{1}|,|\sigma_{2}|\right\}
\geq C|k_{min}||k_{max}|^{2j},\nonumber\\
&&(c):\quad |\sigma_{2}|={\rm max}\left\{|\sigma|,|\sigma_{1}|,|\sigma_{2}|\right\}
\geq C|k_{min}||k_{max}|^{2j}.\nonumber
\end{eqnarray}
When $(a):\quad |\sigma|={\rm max}\left\{|\sigma|,|\sigma_{1}|,|\sigma_{2}|\right\}
\geq C|k_{min}||k_{max}|^{2j}$, we have that
\begin{eqnarray}
K_{1}(k_{1},\tau_{1},k,\tau)=\frac{|kk_{1}k_{2}|\langle k\rangle^{s}}{(1+k^{2})
\langle\sigma\rangle^{1/2}\prod_{l=1}^{2}
\langle k_{l}\rangle^{s}\langle\sigma_{l}\rangle^{1/2}}\leq C
\frac{|k|^{s-\frac{3}{2}}\prod_{l=1}^{2}
k_{l}^{\frac{2-j}{2}-s}}{\prod_{l=1}^{2}\langle\sigma_{l}\rangle^{1/2}}\label{3.04};
\end{eqnarray}
if $\frac{2-j}{2}\leq s\leq \frac{3}{2},$ from (\ref{3.04}), we have that
\begin{eqnarray}
K_{1}(k_{1},\tau_{1},k,\tau)\leq\frac{ C}{\prod_{l=1}^{2}\langle\sigma_{l}\rangle^{1/2}}\label{3.05};
\end{eqnarray}
if $s\geq \frac{3}{2},$ since $s\geq \frac{2-j}{2}$, we have that
\begin{eqnarray}
&&K_{1}(k_{1},\tau_{1},k,\tau)\leq C\frac{|k|^{s-\frac{3}{2}}\left[{\rm max} \left\{|k_{1}|,|k_{2}|\right\}\right]^{\frac{2-j}{2}-s}\left[{\rm min} \left\{|k_{1}|,|k_{2}|\right\}\right]^{\frac{2-j}{2}-s}}{\prod_{l=1}^{2}\langle\sigma_{l}\rangle^{1/2}}\nonumber\\
&&\leq C\frac{\left[{\rm max} \left\{|k_{1}|,|k_{2}|\right\}\right]^{-\frac{1+j}{2}}\left[{\rm min}
 \left\{|k_{1}|,|k_{2}|\right\}\right]^{\frac{2-j}{2}-s}}{\prod_{l=1}^{2}\langle\sigma_{l}\rangle^{1/2}}\nonumber\\&&
\leq \frac{C}{\prod_{l=1}^{2}\langle\sigma_{l}\rangle^{1/2}}\label{3.06};
\end{eqnarray}
from (\ref{3.05})-(\ref{3.06}),
by using  the Plancherel identity and the H\"older inequality
as well as  Lemma 2.11,
we have that
\begin{eqnarray}
&&\int_{\SR^{2}_{\tau k}}\int_{\!\!\!\mbox{\scriptsize $
\begin{array}{l}
k=k_{1}+k_{2}\\
\tau=\tau_{1}+\tau_{2}
\end{array}
$}}K_{1}(k_{1},\tau_{1},k,\tau)F(k,\tau)\prod_{l=1}^{2}F_{l}(k_{l},\tau_{l})
dk_{1}d\tau_{1}dkd\tau \nonumber\\&&
\leq C\int_{\!\!\!\mbox{\scriptsize $
\begin{array}{l}
k=k_{1}+k_{2}\\
\tau=\tau_{1}+\tau_{2}
\end{array}
$}}\frac{F(k,\tau)\prod_{l=1}^{2}F_{l}(k_{l},\tau_{l})}{
\prod_{l=1}^{2}\langle\sigma_{l}\rangle^{1/2}}dk_{1}d\tau_{1}dkd\tau  \nonumber\\&&
\leq C\left\|\mathscr{F}^{-1}\left(F\right)\right\|_{L^{2}_{xt}}\prod_{l=1}^{2}
\left\|\mathscr{F}^{-1}\left(\frac{F_{l}}{\langle\sigma_{l}\rangle^{1/2}}\right)\right\|_{L^{4}_{xt}}\nonumber\\
&&\leq C\delta^{\frac{j}{2j+1}-2\epsilon}\|F\|_{L^{2}}\prod_{l=1}^{2}\|F_{l}\|_{L^{2}}\label{3.07}.
\end{eqnarray}
When $(b):\quad |\sigma_{1}|={\rm max}\left\{|\sigma|,|\sigma_{1}|,|\sigma_{2}|\right\}
\geq C|k_{min}||k_{max}|^{2j}$,
by using the proof similar to (\ref{3.05})-(\ref{3.06}), we have that
\begin{eqnarray}
&&K_{1}(k_{1},\tau_{1},k,\tau)\leq \frac{C}{\langle\sigma\rangle^{1/2}\langle\sigma_{2}\rangle^{1/2}}\label{3.08};
\end{eqnarray}
by using the Cauchy-Schwarz inequality and Lemma 2.12,   we have that
\begin{eqnarray*}
&&\int_{\SR^{2}_{\tau k}}F(k,\tau)
\left(\langle\sigma\rangle^{-1/2}\int_{\!\!\!\mbox{\scriptsize $
\begin{array}{l}
k=k_{1}+k_{2}\\
\tau=\tau_{1}+\tau_{2}
\end{array}
$}}\frac{\prod_{l=1}^{2}F_{l}(k_{l},\tau_{l})}{\langle\sigma_{2}\rangle^{1/2}}dk_{1}
d\tau_{1} \right)dkd\tau \nonumber\\
&&\int_{\SR^{2}_{\tau k}}F(k,\tau)
\left(\langle\sigma\rangle^{-\frac{1}{2}+\epsilon}\int_{\!\!\!\mbox{\scriptsize $
\begin{array}{l}
k=k_{1}+k_{2}\\
\tau=\tau_{1}+\tau_{2}
\end{array}
$}}\frac{\prod_{l=1}^{2}F_{l}(k_{l},\tau_{l})}{\langle\sigma_{2}\rangle^{1/2}}dk_{1}
d\tau_{1} \right)dkd\tau \nonumber\\
&&\leq C
\left\|F(k,\tau)\right\|_{L^{2}_{k\tau}}
\left\|\langle\sigma\rangle^{-1/2+\epsilon}\int_{\!\!\!\mbox{\scriptsize $
\begin{array}{l}
k=k_{1}+k_{2}\\
\tau=\tau_{1}+\tau_{2}
\end{array}
$}}\frac{\prod_{l=1}^{2}F_{l}(k_{l},\tau_{l})}{\langle\sigma_{2}\rangle^{1/2}}dk_{1}
d\tau_{1} \right\|_{L_{k\tau}^{2}}\nonumber\\
&&\leq C\delta^{\frac{j}{2j+1}-2\epsilon}\|F\|_{L^{2}}\prod_{l=1}^{2}\|F_{l}\|_{L^{2}}.
\end{eqnarray*}
When $(c):\quad |\sigma_{2}|={\rm max}\left\{|\sigma|,|\sigma_{1}|,
|\sigma_{2}|\right\}
\geq C|k_{min}||k_{max}|^{2j}$, this case can be proved similarly to case
$(b):\quad |\sigma_{1}|={\rm max}\left\{|\sigma|,|\sigma_{1}|,|\sigma_{2}|\right\}
\geq C|k_{min}||k_{max}|^{2j}$.

We have completed the proof of Lemma 3.1.

\begin{Lemma}\label{Lemma3.2}
Let $u_{l}(x,t)$ with $l=1,2$ which are zero $x$-mean for all $t$ be
$2\pi $- periodic functions  of $x$
and $s\geq -\frac{j}{2}$. For $\epsilon<\frac{1}{10000(2j+1)},$
then we have that
\begin{eqnarray}
      \left\|\partial_{x}\left[
      \prod_{l=1}^{2}\left[\eta\left(\frac{t}{\delta}\right)u_{l}\right]\right]\right\|_{X_{s,-\frac{1}{2}}^{\delta}}\leq
      C\delta
      ^{\frac{j}{2j+1}-2\epsilon}
      \prod_{l=1}^{2}\|u_{l}\|_{X_{s,\frac{1}{2}}^{\delta}}.
        \label{3.016}
\end{eqnarray}
\end{Lemma}
{\bf Proof.}Let $\tilde{u}$ and $\tilde{u}_{1},\tilde{u}_{2}$ be
the extension of $u,u_{1},u_{2}$,
 respectively, according to Lemma 2.10, we have that
\begin{eqnarray*}
\|u\|_{X_{s,\frac{1}{2}}^{\delta}}=\|\tilde{u}\|_{X_{s,\frac{1}{2}}},
\|u_{l}\|_{X_{s,\frac{1}{2}}^{\delta}}
=\|\tilde{u}_{l}\|_{X_{s,\frac{1}{2}}},\quad l=1,2.
\end{eqnarray*}
By duality and  the Plancherel identity,
for $u\in X_{-s,\frac{1}{2}}^{\delta},$
it suffices to prove that
\begin{eqnarray}
&&\int_{\!\!\!\mbox{\scriptsize $
\begin{array}{l}
k=k_{1}+k_{2}\\
\tau=\tau_{1}+\tau_{2}
\end{array}
$}}\left|k\mathscr{F}\left(\eta\left(\frac{t}{\delta}\right)\tilde{u}\right)
(k,\tau)\prod_{l=1}^{2}
\mathscr{F}\left(\eta\left(\frac{t}{\delta}\right)\tilde{u}_{l}\right)(k_{l},\tau_{l})
\right|dk_{1}
d\tau_{1}dkd\tau \nonumber\\&&\leq C\delta
      ^{\frac{j}{2j+1}-2\epsilon}
\|u\|_{X_{-s,\frac{1}{2}}^{\delta}}
\prod_{l=1}^{2}\|u_{l}\|_{X_{s,\frac{1}{2}}^{\delta}}=C\delta
      ^{\frac{j}{2j+1}-2\epsilon}\|\tilde{u}\|_{X_{-s,\frac{1}{2}}}\prod_{l=1}^{2}
      \|\tilde{u}_{l}\|_{X_{s,\frac{1}{2}}}.\label{3.017}
\end{eqnarray}
Without loss of generality, we can assume that
 $\mathscr{F}\left(\eta\left(\frac{t}{\delta}\right)\tilde{u}_{l}\right)
 (k_{l},\tau_{l})\geq 0(l=1,2)$
 and $\mathscr{F}\left(\eta\left(\frac{t}{\delta}\right)\tilde{u}\right)
 (k,\tau)\geq 0$.
Let
\begin{eqnarray*}
&&F(k,\tau)=\langle k\rangle^{-s}\langle\sigma\rangle^{1/2}\mathscr{F}
\left(\eta\left(\frac{t}{\delta}\right)\tilde{u}\right)(k,\tau),\\
&&F_{l}(k_{l},\tau_{l})=\langle k_{l}\rangle^{s}\langle\sigma_{l}\rangle^{1/2}
\mathscr{F}\left(\eta\left(\frac{t}{\delta}\right)\tilde{u}_{l}\right)(k_{l},\tau_{l})
,\quad l=1,2,\\
&&K_{2}(k_{1},\tau_{1},k,\tau)=\frac{|k|\langle k\rangle^{s}}{
\langle\sigma\rangle^{1/2}
\prod_{l=1}^{2}\langle k_{l}\rangle^{s}\langle\sigma_{l}\rangle^{1/2}}.\\
\end{eqnarray*}
To obtain (\ref{3.017}),  it suffices to prove that
\begin{eqnarray}
&&\int_{\!\!\!\mbox{\scriptsize $
\begin{array}{l}
k=k_{1}+k_{2}\\
\tau=\tau_{1}+\tau_{2}
\end{array}
$}}K_{2}(k_{1},\tau_{1},k,\tau)F(k,\tau)\prod_{l=1}^{2}F_{l}(k_{l},\tau_{l})
dk_{1}d\tau_{1} dkd\tau\nonumber
\\&&\leq C\delta
      ^{\frac{j}{2j+1}-2\epsilon}\|F\|_{L^{2}}\prod_{l=1}^{2}\|F_{l}\|_{L^{2}}\label{3.018}.
\end{eqnarray}
From the mean zero condition, we can assume that $k\neq0,k_{l}\neq0(l=1,2).$
Since ${\rm min}\left\{|k|,|k_{1}|,|k_{2}|\right\}\geq 1,$ from Lemma 2.4,
we have that one of the following three cases
\begin{eqnarray}
&&(a):\quad |\sigma|={\rm max}\left\{|\sigma|,|\sigma_{1}|,|\sigma_{2}|\right\}
\geq C|k_{min}||k_{max}|^{2j},\nonumber\\
&&(b):\quad |\sigma_{1}|={\rm max}\left\{|\sigma|,|\sigma_{1}|,|\sigma_{2}|\right\}
\geq C|k_{min}||k_{max}|^{2j},\nonumber\\
&&(c):\quad |\sigma_{2}|={\rm max}\left\{|\sigma|,|\sigma_{1}|,|\sigma_{2}|\right\}
\geq C|k_{min}||k_{max}|^{2j}.\nonumber
\end{eqnarray}
When $(a):\quad |\sigma|={\rm max}\left\{|\sigma|,|\sigma_{1}|,|\sigma_{2}|\right\}
\geq C|k_{min}||k_{max}|^{2j}$,  we have that
\begin{eqnarray}
K_{2}(k_{1},\tau_{1},k,\tau)=\frac{|k|\langle k\rangle^{s}}{
\langle\sigma\rangle^{1/2}\prod_{l=1}^{2}
\langle k_{l}\rangle^{s}\langle\sigma_{l}\rangle^{1/2}}\leq C
\frac{|k|^{s+\frac{1}{2}}\prod_{l=1}^{2}
k_{l}^{-\frac{j}{2}-s}}{\prod_{l=1}^{2}\langle\sigma_{l}\rangle^{1/2}}\label{3.019};
\end{eqnarray}
if $-\frac{j}{2}\leq s\leq -\frac{1}{2},$ from (\ref{3.019}), we have that
\begin{eqnarray}
K_{2}(k_{1},\tau_{1},k,\tau)\leq\frac{ C}{\prod_{l=1}^{2}\langle\sigma_{l}\rangle^{1/2}}\label{3.020};
\end{eqnarray}
if $s\geq -\frac{1}{2},$ since $s\geq -\frac{j}{2}$, we have that
\begin{eqnarray}
&&K_{2}(k_{1},\tau_{1},k,\tau)\leq C\frac{|k|^{s+\frac{1}{2}}\left[{\rm max} \left\{|k_{1}|,|k_{2}|\right\}\right]^{-\frac{j}{2}-s}\left[{\rm min} \left\{|k_{1}|,|k_{2}|\right\}\right]^{-\frac{j}{2}-s}}{\prod_{l=1}^{2}\langle\sigma_{l}\rangle^{1/2}}\nonumber\\
&&\leq C\frac{\left[{\rm max} \left\{|k_{1}|,|k_{2}|\right\}\right]^{-\frac{j-1}{2}}\left[{\rm min} \left\{|k_{1}|,|k_{2}|\right\}\right]^{-\frac{j}{2}-s}}{\prod_{l=1}^{2}\langle\sigma_{l}\rangle^{1/2}}\nonumber\\&&
\leq \frac{C}{\prod_{l=1}^{2}\langle\sigma_{l}\rangle^{1/2}}\label{3.021};
\end{eqnarray}
from (\ref{3.020})-(\ref{3.021}),
by using  the Plancherel identity and the H\"older inequality
and Lemma 2.11,
we have that
\begin{eqnarray}
&&\int_{\SR^{2}_{\tau k}}\int_{\!\!\!\mbox{\scriptsize $
\begin{array}{l}
k=k_{1}+k_{2}\\
\tau=\tau_{1}+\tau_{2}
\end{array}
$}}K_{2}(k_{1},\tau_{1},k,\tau)F(k,\tau)\prod_{l=1}^{2}F_{l}(k_{l},\tau_{l})
dk_{1}d\tau_{1}dkd\tau \nonumber\\&&
\leq C\int_{\!\!\!\mbox{\scriptsize $
\begin{array}{l}
k=k_{1}+k_{2}\\
\tau=\tau_{1}+\tau_{2}
\end{array}
$}}\frac{F(k,\tau)\prod_{l=1}^{2}F_{l}(k_{l},\tau_{l})}{
\prod_{l=1}^{2}\langle\sigma_{l}\rangle^{1/2}}dk_{1}d\tau_{1}dkd\tau  \nonumber\\&&
\leq C\left\|\mathscr{F}^{-1}\left(F\right)\right\|_{L^{2}_{xt}}\prod_{l=1}^{2}
\left\|\mathscr{F}^{-1}\left(\frac{F_{l}}{\langle\sigma_{l}\rangle^{1/2}}\right)\right\|_{L^{4}_{xt}}\nonumber\\
&&\leq C\delta^{\frac{j}{2j+1}-2\epsilon}\|F\|_{L^{2}}\prod_{l=1}^{2}\|F_{l}\|_{L^{2}}\label{3.022}.
\end{eqnarray}
When $(b):\quad |\sigma_{1}|={\rm max}\left\{|\sigma|,|\sigma_{1}|,|\sigma_{2}|\right\}
\geq C|k_{min}||k_{max}|^{2j}$,
by using the proof similar to (\ref{3.020})-(\ref{3.021}), we have that
\begin{eqnarray}
&&K_{2}(k_{1},\tau_{1},k,\tau)\leq \frac{C}{\langle\sigma\rangle^{1/2}\langle\sigma_{2}\rangle^{1/2}}\label{3.023};
\end{eqnarray}
by using the Cauchy-Schwarz inequality and Lemma 2.12,   we have that
\begin{eqnarray}
&&\int_{\SR^{2}_{\tau k}}\int_{\!\!\!\mbox{\scriptsize $
\begin{array}{l}
k=k_{1}+k_{2}\\
\tau=\tau_{1}+\tau_{2}
\end{array}
$}}K_{2}(k_{1},\tau_{1},k,\tau)F(k,\tau)\prod_{l=1}^{2}F_{l}(k_{l},\tau_{l})
dk_{1}d\tau_{1}dkd\tau \nonumber\\
&&\leq C\int_{\SR^{2}_{\tau k}}F(k,\tau)
\left(\langle\sigma\rangle^{-1/2}\int_{\!\!\!\mbox{\scriptsize $
\begin{array}{l}
k=k_{1}+k_{2}\\
\tau=\tau_{1}+\tau_{2}
\end{array}
$}}\frac{\prod_{l=1}^{2}F_{l}(k_{l},\tau_{l})}{\langle\sigma_{2}\rangle^{1/2}}dk_{1}
d\tau_{1} \right)dkd\tau \nonumber\\
&&\leq C\int_{\SR^{2}_{\tau k}}F(k,\tau)
\left(\langle\sigma\rangle^{-\frac{1}{2}+\epsilon}\int_{\!\!\!\mbox{\scriptsize $
\begin{array}{l}
k=k_{1}+k_{2}\\
\tau=\tau_{1}+\tau_{2}
\end{array}
$}}\frac{\prod_{l=1}^{2}F_{l}(k_{l},\tau_{l})}{\langle\sigma_{2}\rangle^{1/2}}dk_{1}
d\tau_{1} \right)dkd\tau \nonumber\\
&&\leq C\|F\|_{L^{2}}\left\|
\langle\sigma\rangle^{-\frac{1}{2}+\epsilon}\int_{\!\!\!\mbox{\scriptsize $
\begin{array}{l}
k=k_{1}+k_{2}\\
\tau=\tau_{1}+\tau_{2}
\end{array}
$}}\frac{\prod_{l=1}^{2}F_{l}(k_{l},\tau_{l})}{\langle\sigma_{2}\rangle^{1/2}}dk_{1}
d\tau_{1} \right\|_{L^{2}}\nonumber\\
&&\leq C\delta^{\frac{j}{2j+1}-2\epsilon}\|F\|_{L^{2}}\prod_{l=1}^{2}\|F_{l}\|_{L^{2}}.\nonumber
\end{eqnarray}
When $(c):\quad |\sigma_{2}|={\rm max}\left\{|\sigma|,|\sigma_{1}|,
|\sigma_{2}|\right\}
\geq C|k_{min}||k_{max}|^{2j}$, this case can be proved similarly to case
$(b):\quad |\sigma_{1}|={\rm max}\left\{|\sigma|,|\sigma_{1}|,|\sigma_{2}|\right\}
\geq C|k_{min}||k_{max}|^{2j}$.

We have completed the proof of Lemma 3.2.

\begin{Lemma}\label{Lemma3.3}
Let $u_{l}(x,t)$ with $l=1,2$ which are zero $x$-mean for all $t$ be
$2\pi$-periodic functions  of $x$
and $s\geq -\frac{j+2}{2}$. For $\epsilon<\frac{1}{10000(2j+1)},$
then we have that
\begin{eqnarray}
      \left\|\partial_{x}(1-\partial_{x}^{2})^{-1}\left[
      \prod_{l=1}^{2}\left[\eta\left(\frac{t}{\delta}\right)u_{l}\right]\right]\right\|_{X_{s,-\frac{1}{2}}^{\delta}}\leq
      C\delta
      ^{\frac{j}{2j+1}-2\epsilon}
      \prod_{l=1}^{2}\|u_{l}\|_{X_{s,\frac{1}{2}}^{\delta}}.
        \label{3.031}
\end{eqnarray}
\end{Lemma}

Lemma 3.3 can be proved similarly to Lemma 3.2.

\begin{Lemma}\label{Lemma3.4}
Let $v_{l}(x,t)$ with $l=1,2$ which are zero $x$-mean for all $t$ be
$2\pi$-periodic functions  of $x$. For $s\geq \frac{2-j}{2}$, we have that
\begin{eqnarray}
      &&\left\|\int_{\!\!\!\mbox{\scriptsize $
\begin{array}{l}
k=k_{1}+k_{2}\\
\tau=\tau_{1}+\tau_{2}
\end{array}
$}}\frac{|kk_{1}k_{2}|\langle k\rangle ^{s}}{\langle\sigma\rangle(1+k^{2})}\prod_{l=1}^{2}\mathscr{F}\left(\eta\left(\frac{t}{\delta}\right)\tilde{v}_{l}\right)(k_{l},\tau_{l})dk_{1}d\tau_{1}
      \right\|_{(L^{2}(k)L^{1}(\tau))}\nonumber\\&&\leq C\delta
      ^{\frac{j}{2j+1}-3\epsilon}\prod_{l=1}^{2}\|v_{l}\|_{X_{s,\frac{1}{2}}^{\delta}}.
        \label{3.032}
\end{eqnarray}
\end{Lemma}
{\bf Proof.}
Let $\tilde{v}_{1},\tilde{v}_{2}$ be
the extension of $v_{1},v_{2}$,
 respectively, according to Lemma 2.10, we have that
\begin{eqnarray*}
\|v_{l}\|_{X_{s,\frac{1}{2}}^{\delta}}
=\|\tilde{v}_{l}\|_{X_{s,\frac{1}{2}}},\quad l=1,2.
\end{eqnarray*}
Without loss of generality, we can assume that
 $\mathscr{F}\left(\eta\left(\frac{t}{\delta}\right)\tilde{v}_{l}\right)(k_{l},\tau_{l})\geq 0(l=1,2)$
 and $\mathscr{F}\left(\eta\left(\frac{t}{\delta}\right)\tilde{v}\right)(k,\tau)\geq 0$.
Let
\begin{eqnarray*}
&&G_{l}(k_{l},\tau_{l})=\langle k_{l}\rangle^{s}\langle\sigma_{l}\rangle^{1/2}\mathscr{F}
\left(\eta\left(\frac{t}{\delta}\right)\tilde{v}_{l}\right)(k_{l},\tau_{l})
,\quad l=1,2,\\
&&K_{3}(k_{1},\tau_{1},k,\tau)=\frac{|kk_{1}k_{2}|\langle k\rangle^{s}}{(1+k^{2})
\langle\sigma\rangle
\prod_{l=1}^{2}\langle k_{l}\rangle^{s}\langle\sigma_{l}\rangle^{1/2}}.\\
\end{eqnarray*}
To obtain (\ref{3.032}), it suffices to prove that
\begin{eqnarray}
     && \left\|\int_{\!\!\!\mbox{\scriptsize $
\begin{array}{l}
k=k_{1}+k_{2}\\
\tau=\tau_{1}+\tau_{2}
\end{array}
$}}K_{3}(k_{1},\tau_{1},k,\tau)
     \prod_{l=1}^{2}G_{l}(k_{l},\tau_{l})dk_{1}d\tau_{1}
      \right\|_{(L^{2}(k)L^{1}(\tau))}\nonumber\\&&\leq C\delta
      ^{\frac{j}{2j+1}-3\epsilon}\prod_{l=1}^{2}\|G_{l}\|_{L^{2}}.
        \label{3.033}
\end{eqnarray}
Since  ${\rm min}\left\{|k|,|k_{1}|,|k_{2}|\right\}\geq 1,$
from Lemma 2.4, we know that one of the following three  cases must occur:
\begin{eqnarray*}
&&(a):\quad |\sigma|={\rm max}\left\{|\sigma|,|\sigma_{1}|,|\sigma_{2}|\right\}
\geq C|k_{min}||k_{max}|^{2j},\nonumber\\
&&(b):\quad |\sigma_{1}|={\rm max}\left\{|\sigma|,|\sigma_{1}|,|\sigma_{2}|\right\}
\geq C|k_{min}||k_{max}|^{2j},\nonumber\\
&&(c):\quad |\sigma_{2}|={\rm max}\left\{|\sigma|,|\sigma_{1}|,|\sigma_{2}|\right\}
\geq C|k_{min}||k_{max}|^{2j}.\nonumber
\end{eqnarray*}

\noindent When $(a):\quad |\sigma|={\rm max}\left\{|\sigma|,|\sigma_{1}|,|\sigma_{2}|\right\}
\geq C|k_{min}||k_{max}|^{2j}$. If $\langle\sigma_{1}\rangle\geq
C|k_{min}|^{\epsilon^{'}}|k_{max}|^{2j\epsilon^{'}}$,
in this case,
by using the proof similar to  (\ref{3.05})-(\ref{3.06}), we have that
\begin{eqnarray}
K_{3}(k_{1},\tau_{1},k,\tau)\leq \frac{|k|^{s-\frac{3}{2}}\prod_{l=1}^{2}\langle k_{l}\rangle^{\frac{2-j}{2}-s}}
{\langle \sigma \rangle^{\frac{1}{2}+\epsilon\epsilon^{'}}\langle \sigma_{1}
 \rangle^{\frac{1}{2}-\epsilon}\langle \sigma_{2} \rangle^{\frac{1}{2}}}\leq \frac{C}{\langle \sigma
 \rangle^{\frac{1}{2}+\epsilon\epsilon^{'}}\langle \sigma_{1} \rangle^{\frac{1}{2}-\epsilon}\langle \sigma_{2}
 \rangle^{\frac{1}{2}}};\label{3.034}
\end{eqnarray}
by using (\ref{3.034}),  the Cauchy-Schwarz inequality and the Plancherel identity as well as
Lemmas 2.3, 2.13, then we have that
\begin{eqnarray*}
      &&\left\|\langle\sigma\rangle^{-\frac{1}{2}-\epsilon\epsilon^{'}}\int_{\!\!\!\mbox{\scriptsize $
\begin{array}{l}
k=k_{1}+k_{2}\\
\tau=\tau_{1}+\tau_{2}
\end{array}
$}}\frac{\prod_{l=1}^{2}G_{l}(k_{l},\tau_{l})}{\langle \sigma_{1}
 \rangle^{\frac{1}{2}-\epsilon}\langle\sigma_{2}\rangle^{1/2}}dk_{1}d\tau_{1}
      \right\|_{L^{2}(k)L^{1}(d\tau)}\nonumber\\&&\leq C
      \left\|\int_{\!\!\!\mbox{\scriptsize $
\begin{array}{l}
k=k_{1}+k_{2}\\
\tau=\tau_{1}+\tau_{2}
\end{array}
$}}\frac{\prod_{l=1}^{2}G_{l}(k_{l},\tau_{l})}{\langle \sigma_{1} \rangle^{\frac{1}{2}-\epsilon}\langle\sigma_{2}
\rangle^{1/2}}dk_{1}d\tau_{1}
\right\|_{L_{k}^{2}L_{\tau}^{2}}\nonumber\\
&&\leq C\left\|\mathscr{F}^{-1}
\left(\frac{G_{1}}{\langle\sigma_{1}\rangle^{\frac{1}{2}-\epsilon}}\right)\right\|_{L_{xt}^{4}}\left\|\mathscr{F}^{-1}
\left(\frac{G_{2}}{\langle\sigma_{2}\rangle^{\frac{1}{2}}}\right)\right\|_{L_{xt}^{4}}\\
&&\leq C\left\|\eta\left(\frac{t}{\delta}\right)\tilde{v}_{1}\right\|_{X_{s,\frac{j+1}{2(2j+1)}+\epsilon}}
\left\|\eta\left(\frac{t}{\delta}\right)\tilde{v}_{2}\right\|_{X_{s,\frac{j+1}{2(2j+1)}}}\nonumber\\
&&\leq C\delta^{\frac{j}{2(2j+1)}-3\epsilon}\prod_{l=1}^{2}\left\|\tilde{v}_{l}\right\|_{X_{s,\frac{1}{2}-\epsilon}}\nonumber\\
&&\leq C\delta^{\frac{j}{2(2j+1)}-3\epsilon}\prod_{l=1}^{2}\left\|\tilde{v}_{l}\right\|_{X_{s,\frac{1}{2}}}\nonumber\\
&&\leq C\delta^{\frac{j}{2j+1}-3\epsilon}
       \prod_{l=1}^{2}\|G_{l}\|_{L^{2}};
        \label{3.035}
\end{eqnarray*}
If $\langle\sigma_{2}\rangle\geq C|k_{min}|^{\epsilon^{'}}|k_{max}|^{2j\epsilon^{'}}$, this case can be proved similarly to case
 $\langle\sigma_{1}\rangle\geq C|k_{min}|^{\epsilon}|k_{max}|^{\epsilon^{'}}$.
if $\langle\sigma_{l}\rangle\leq C|k_{min}|^{\epsilon}|k_{max}|^{2j\epsilon^{'}},l=1,2,$ in this case we have that
\begin{eqnarray}
\mu=k^{2j+1}-k_{1}^{2j+1}-k_{2}^{2j+1}
+O(\langle |k_{min}||k_{max}|^{2j}\rangle^{\epsilon^{'}})\label{3.036}
\end{eqnarray}
and
\begin{eqnarray*}
K_{3}(k_{1},\tau_{1},k,\tau)\leq C\frac{|k|^{s-\frac{3}{2}}
\prod_{l=1}^{2}|k_{l}|^{\frac{2-j}{2}-s}}{\langle\sigma\rangle^{1/2}\prod_{l=1}^{2}\langle\sigma_{l}\rangle^{1/2}}
\end{eqnarray*}
by using the proof similar to (\ref{3.05})-(\ref{3.06}), we have that
\begin{eqnarray}
K_{3}(k_{1},\tau_{1},k,\tau)\leq \frac{C}{\langle\sigma\rangle^{1/2}\prod_{l=1}^{2}\langle\sigma_{l}\rangle^{1/2}}\label{3.037}.
\end{eqnarray}
Consequently, by using (\ref{3.021}) and the Cauchy-Schwartz inequality with respect to $\tau$ and Lemmas 2.8,  2.11, we have that
\begin{eqnarray}
    &&  \left\|\langle\sigma\rangle^{-\frac{1}{2}}\chi_{\Omega(k)}\int_{\!\!\!\mbox{\scriptsize $
\begin{array}{l}
k=k_{1}+k_{2}\\
\tau=\tau_{1}+\tau_{2}
\end{array}
$}}\frac{\prod_{l=1}^{2}G_{l}(k_{l},\tau_{l})}{\prod_{l=1}^{2}\langle\sigma_{l}
\rangle^{1/2}}dk_{1}d\tau_{1}
      \right\|_{L^{2}(kL^{1}(d\tau)}\nonumber\\&&\leq C\left\|\left(\int\langle\sigma\rangle^{-1}
      \chi_{\Omega(k)}(\mu)d\tau\right)^{1/2}\int_{\!\!\!\mbox{\scriptsize $
\begin{array}{l}
k=k_{1}+k_{2}\\
\tau=\tau_{1}+\tau_{2}
\end{array}
$}}
      \frac{\prod_{l=1}^{2}G_{l}(k_{l},\tau_{l})}{\prod_{l=1}^{2}\langle\sigma_{l}\rangle^{1/2}}dk_{1}d\tau_{1}\right\|
      _{L_{ k\tau}^{2}}\nonumber\\
      &&\leq C\left(\int\langle\sigma\rangle^{-1}
      \chi_{\Omega(k)}(\mu)d\tau\right)^{1/2}\left\|\int_{\!\!\!\mbox{\scriptsize $
\begin{array}{l}
k=k_{1}+k_{2}\\
\tau=\tau_{1}+\tau_{2}
\end{array}
$}}
      \frac{\prod_{l=1}^{2}G_{l}(k_{l},\tau_{l})}{\prod_{l=1}^{2}\langle\sigma_{l}\rangle^{1/2}}dk_{1}d\tau_{1}\right\|
      _{L_{ k\tau}^{2}}\nonumber\\
      &&\leq C\delta^{\frac{j}{2j+1}-2\epsilon}\prod_{l=1}^{2}\|G_{l}\|_{L^{2}},
        \label{3.038}
\end{eqnarray}
where
\begin{eqnarray*}
\Omega(k)=\left\{\mu \in \R: \quad |\mu|\sim M, \mu=C|k_{min}||k_{max}|^{2j}
+O(\langle |k_{min}||k_{max}|^{2j}\rangle^{\epsilon^{'}})\right\}.
\end{eqnarray*}
When $(b):\quad |\sigma_{1}|={\rm max}\left\{|\sigma|,|\sigma_{1}|,|\sigma_{2}|\right\}
\geq C|k_{min}||k_{max}|^{2j}$.
by using the proof similar to  (\ref{3.05})-(\ref{3.06}), we have that
\begin{eqnarray}
K_{3}(k_{1},\tau_{1},k,\tau)\leq C\frac{|k|^{s-\frac{3}{2}}\prod_{l=1}^{2}\langle k_{l}\rangle^{\frac{2-j}{2}-s}}
{\langle \sigma \rangle \langle \sigma_{2} \rangle^{\frac{1}{2}}}\leq \frac{C}{\langle \sigma
 \rangle\langle \sigma_{2}
 \rangle^{\frac{1}{2}}},\label{3.039}
\end{eqnarray}
by using the Cauchy-Schwarz inequality with respect to $\tau$ and Lemma 2.12,  we have that
\begin{eqnarray*}
  && \left\|\int_{\!\!\!\mbox{\scriptsize $
\begin{array}{l}
k=k_{1}+k_{2}\\
\tau=\tau_{1}+\tau_{2}
\end{array}
$}}\frac{\prod_{l=1}^{2}G_{l}(k_{l},\tau_{l})}{\langle\sigma\rangle\langle\sigma_{2}\rangle^{1/2}}dk_{1}d\tau_{1}
      \right\|_{L^{2}(k)L^{1}(\tau)}\nonumber\\&&\leq C\left\|\langle \sigma\rangle^{-\frac{1}{2}+\epsilon}\int_{\!\!\!\mbox{\scriptsize $
\begin{array}{l}
k=k_{1}+k_{2}\\
\tau=\tau_{1}+\tau_{2}
\end{array}
$}}\frac{\prod_{l=1}^{2}G_{l}(k_{l},\tau_{l})}{\langle\sigma_{2}\rangle^{1/2}}dk_{1}
d\tau_{1}\right\|_{L_{k\tau}^{2}}
\nonumber\\
&&\leq C\delta^{\frac{j}{2j+1}-2\epsilon}
\prod_{l=1}^{2}\|G_{l}\|_{L^{2}}.
\end{eqnarray*}
When $(c):\quad |\sigma_{2}|={\rm max}\left\{|\sigma|,|\sigma_{1}|,
|\sigma_{2}|\right\}
\geq C|k_{min}||k_{max}|^{2j}$. This case can be proved similarly to  $(b):\quad
|\sigma_{1}|={\rm max}\left\{|\sigma|,|\sigma_{1}|,|\sigma_{2}|\right\}
\geq C|k_{min}||k_{max}|^{2j}$.

We have completed the proof Lemma 3.4.

\begin{Lemma}\label{Lemma3.5}
Let $v_{l}(x,t)$ with $l=1,2$ which are zero $x$-mean for all $t$ be
$2\pi$- periodic functions  of $x$. For $s\geq -\frac{j}{2}$, we have that
\begin{eqnarray}
      &&\left\|\int_{\!\!\!\mbox{\scriptsize $
\begin{array}{l}
k=k_{1}+k_{2}\\
\tau=\tau_{1}+\tau_{2}
\end{array}
$}}\frac{|k|\langle k\rangle ^{s}}{\langle\sigma\rangle}\prod_{l=1}^{2}\mathscr{F}
\left(\eta\left(\frac{t}{\delta}\right)v_{l}\right)(k_{l},\tau_{l})dk_{1}d\tau_{1}
      \right\|_{(L^{2}(k)L^{1}(\tau))}\nonumber\\&&\leq C\delta
      ^{\frac{j}{2j+1}-2\epsilon}\prod_{l=1}^{2}\|v_{l}\|_{X_{s,\frac{1}{2}}^{\delta}}.
        \label{3.044}
\end{eqnarray}
\end{Lemma}

By using the proof similar to Lemma 3.4, we can obtain Lemma 3.5.

\begin{Lemma}\label{Lemma3.6}
Let $v_{l}(x,t)$ with $l=1,2$ which are zero $x$-mean for all $t$ be
$2\pi $- periodic functions  of $x$. For $s\geq -\frac{j}{2}$, we have that
\begin{eqnarray}
      &&\left\|\int_{\!\!\!\mbox{\scriptsize $
\begin{array}{l}
k=k_{1}+k_{2}\\
\tau=\tau_{1}+\tau_{2}
\end{array}
$}}\frac{|k|\langle k\rangle ^{s}}{(1+k^{2})\langle\sigma\rangle}\prod_{l=1}^{2}
\mathscr{F}\left(\eta\left(\frac{t}{\delta}\right)v_{l}\right)(k_{l},\tau_{l})dk_{1}d\tau_{1}
      \right\|_{L^{2}(k)L^{1}(\tau)}\nonumber\\&&\leq C\delta
      ^{\frac{j}{2j+1}-2\epsilon}\prod_{l=1}^{2}\|v_{l}\|_{X_{s,\frac{1}{2}}^{\delta}}.
        \label{3.045}
\end{eqnarray}
\end{Lemma}

By using the proof similar to Lemma 3.4, we can obtain Lemma 3.6.

\begin{Lemma}\label{Lemma3.7}
Let $u_{l}(x,t)$ with $l=1,2$  which are zero $x$-mean for all $t$ be  $2\pi$-periodic functions  of $x$. Then
\begin{eqnarray}
      \left\|\partial_{x}(1-\partial_{x}^{2})^{-1}\left[
      \prod_{l=1}^{2}\left[\partial_{x}\eta\left(\frac{t}{\delta}\right)u_{l}\right]\right]\right\|_{Z_{s}}\leq C\delta
      ^{\frac{j}{2j+1}-2\epsilon}\prod_{l=1}^{2}\|u_{l}\|_{Y_{s}^{\delta}}.
        \label{3.046}
\end{eqnarray}
\end{Lemma}
Combining Lemma 3.1 with Lemma 3.4, we have Lemma 3.7.

\begin{Lemma}\label{Lemma3.8}
Let $u(x,t)$ with $l=1,2$  which are zero $x$-mean for all $t$ be  $2\pi$-periodic functions  of $x$. Then
\begin{eqnarray}
      \left\|\partial_{x}\left[
      \prod_{l=1}^{2}\left[\eta\left(\frac{t}{\delta}\right)u_{l}\right]\right]\right\|_{Z_{s}}\leq C\delta
      ^{\frac{j}{2j+1}-2\epsilon}\prod_{l=1}^{2}\|u_{l}\|_{Y_{s}^{\delta}}.
        \label{3.047}
\end{eqnarray}
\end{Lemma}

Combining Lemma 3.2 with Lemma 3.5, we have Lemma 3.8.

\begin{Lemma}\label{Lemma3.9}
Let $u(x,t)$ with $l=1,2$  which are zero $x$-mean for all $t$ be  $2\pi$-periodic functions  of $x$. Then
\begin{eqnarray}
      \left\|\partial_{x}(1-\partial_{x}^{2})^{-1}\left[
      \prod_{l=1}^{2}\left[\eta\left(\frac{t}{\delta}\right)u_{l}\right]\right]\right\|_{Z_{s}}\leq C\delta
      ^{\frac{j}{2j+1}-2\epsilon}\prod_{l=1}^{2}\|u_{l}\|_{Y_{s}^{\delta}}.
        \label{3.048}
\end{eqnarray}
\end{Lemma}

Combining Lemma 3.3 with Lemma 3.6, we have Lemma 3.9.

\bigskip
\bigskip

\noindent {\large\bf 4. Proof of Theorem  1.1}

\setcounter{equation}{0}

 \setcounter{Theorem}{0}

\setcounter{Lemma}{0}

\setcounter{section}{4}
Now we are in a position to prove Theorem 1.1.
We define
\begin{eqnarray}
&&\Phi(u)=\eta(t) S(t)\phi-\eta(t) \int_{0}^{t}S(t-t^{'})\eta (t^{'})
A(u)dt^{'},\label{4.01}\nonumber\\
&&B=\left\{u\in Y_{s}^{\delta}: \quad \|u\|_{ Y_{s}^{\delta}}\leq 2C\|\phi\|_{H^{s}
(\mathbf{T})}\right\},
\end{eqnarray}
where
\begin{eqnarray*}
A(u)=\frac{1}{2}\partial_{x}\left[(\eta\left(\frac{t}{\delta}\right)u)^{2}\right]
+(1-\partial_{x}^{2})^{-1}\left[(\eta\left(\frac{t}{\delta}\right)u)^{2}+\frac{1}{2}(\eta\left(\frac{t}{\delta}\right)u_{x})^{2}\right].
\end{eqnarray*}
By using   Lemmas 2.8-2.9, 3.7-3.9,    for sufficiently small $\delta >0$, we have that
\begin{eqnarray*}
\delta
      ^{\frac{j}{2j+1}-3\epsilon}\|\phi\|_{H^{s}(\mathbf{T})}\leq \frac{1}{4},
\end{eqnarray*}
which yields that
\begin{eqnarray}
&&\left\|\Phi(u)\right\|_{Y_{s}}\leq \left\|\eta(t) S(t)\phi\right\|_{Y_{s}^{\delta}}+
\left\|-\frac{1}{2}\eta(t) \int_{0}^{t}S(t-t^{'})\eta (t^{'})
A(u)dt^{'}\right\|_{Y_{s}^{\delta}}\nonumber\\&&\leq C_{1}
\|\phi\|_{H^{s}(\mathbf{T})}+C\left\|\eta\left(\frac{t}{\delta}\right)A(u)\right\|_{Z_{s}}\nonumber\\
&&\leq C\|\phi\|_{H^{s}(\mathbf{T})}+C\delta
      ^{\frac{j}{2j+1}-3\epsilon}\|u\|_{Y_{s}^{\delta}}^{2}\nonumber\\
      &&\leq C\|\phi\|_{H^{s}(\mathbf{T})}+C\delta
      ^{\frac{j}{2j+1}-3\epsilon}\|\phi\|_{H^{s}(\mathbf{T})}^{2}\leq 2C\|\phi\|_{H^{s}
(\mathbf{T})}\label{4.02}
\end{eqnarray}
For $u,v\in B$,  for sufficiently small $\delta >0$,    we have that
\begin{eqnarray}
&&\left\|\Phi(u)-\Phi(v)\right\|_{Y_{s}^{\delta}}\nonumber\\&&\leq C\delta
      ^{\frac{j}{2j+1}-3\epsilon}
\left(\|u\|_{Y_{s}^{\delta}}+\|v\|_{Y_{s}^{\delta}}\right)\|u-v\|_{Y_{s}^{\delta}}\nonumber\\
&&\leq 2C\delta
      ^{\frac{j}{2j+1}-3\epsilon}\|\phi\|_{H^{s}(\mathbf{T})}\|u-v\|_{Y_{s}^{\delta}}\nonumber\\&&\leq
\frac{1}{2}\|u-v\|_{Y_{s}^{\delta}}.\label{4.03}
\end{eqnarray}
From (\ref{4.03}), by using the fixed point Theorem, we have that there exists a $u$ such that $\Phi(u)=u.$
The proof of the remainder of Theorem 1.1 is standard.

We have completed the proof of Theorem 1.1.

\bigskip
\bigskip

\noindent {\large\bf 5. Modified energy}

\setcounter{equation}{0}

 \setcounter{Theorem}{0}

\setcounter{Lemma}{0}

\setcounter{section}{5}
 In this section, we give the almost conserved law which can be used to extend the local solution to the Cauchy problem for (\ref{1.01}) to the global solution
 to the Cauchy problem for (\ref{1.01}).

\begin{Lemma}\label{Lemma5.1}
Let $\frac{2-j}{2}\leq s <1$ and $u$ be the solution to the Cauchy problem for (\ref{1.01}) on $[0,\delta]$. Then
\begin{eqnarray}
&&\left|\int_{0}^{\delta}\int_{\mathbf{T}}\partial_{x}^{3}(I\eta\left(\frac{t}{\delta}\right)u)
\left[I(\eta\left(\frac{t}{\delta}\right)u)^{2}-(\eta\left(\frac{t}{\delta}\right) Iu)^{2}\right]dxdt\right|\nonumber\\&&
\leq C\delta^{\frac{j}{2j+1}-2\epsilon}N^{-j}\|Iu\|_{X_{1,\frac{1}{2}}^{\delta}}^{3}\label{5.01}.
\end{eqnarray}
\end{Lemma}
{\bf Proof.}
To obtain (\ref{5.01}), it suffices to prove that
\begin{eqnarray}
&&\int_{\!\!\!\mbox{\scriptsize $
\begin{array}{l}
k=k_{1}+k_{2}\\
\tau=\tau_{1}+\tau_{2}
\end{array}
$}}\frac{|k|^{3}|\left|m(k)-m(k_{1})m(k_{2})\right|}{\prod_{l=1}^{2}m(k_{l})}\nonumber\\&&
\qquad \times\left|\mathscr{F}(\eta\left(\frac{t}{\delta}\right)\tilde{u})(\tau,k)
\prod_{l=1}^{2}\mathscr{F}(\eta\left(\frac{t}{\delta}\right)\tilde{u}_{l})(\tau_{l},k_{l})\right|
dk_{1}d\tau_{1} dkd\tau\nonumber\\&&
\leq C\delta^{\frac{j}{2j+1}-2\epsilon}N^{-j}\|\tilde{u}\|_{X_{1,\frac{1}{2}}}\prod_{l=1}^{2}
\|\tilde{u}_{l}\|_{X_{1,\frac{1}{2}}}\label{5.02}
\end{eqnarray}
where
\begin{eqnarray*}
\|\tilde{u}\|_{X_{1,\frac{1}{2}}}=\|u\|_{X_{1,\frac{1}{2}}^{\delta}},\quad
\|\tilde{u}_{l}\|_{X_{1,\frac{1}{2}}}=\|u_{l}\|_{X_{1,\frac{1}{2}}^{\delta}},l=1,2.
\end{eqnarray*}
Let
\begin{eqnarray*}
&&H_{l}(k_{l},\tau_{l})=\langle k_{l}\rangle\langle\sigma_{l}\rangle ^{1/2}\mathscr{F}
\left(\eta\left(\frac{t}{\delta}\right)\tilde{u}_{l}\right)(k_{l},\tau_{l}),l=1,2,
\nonumber\\
&&H(k,\tau)=\langle k\rangle\langle\sigma\rangle ^{1/2}\mathscr{F}
\left(\eta\left(\frac{t}{\delta}\right)\tilde{u}\right)(k,\tau).
\end{eqnarray*}
To prove (\ref{5.02}), it suffices to prove
\begin{eqnarray}
&&\int_{\!\!\!\mbox{\scriptsize $
\begin{array}{l}
k=k_{1}+k_{2}\\
\tau=\tau_{1}+\tau_{2}
\end{array}
$}}\frac{\left|m(k)-m(k_{1})m(k_{2})\right||k|^{3}H(k,\tau)\prod_{l=1}^{2}H_{l}(k_{l},\tau_{l})}
{m(k_{1})m(k_{2})\langle\sigma\rangle ^{1/2}\langle k\rangle
\prod_{l=1}^{2}\langle\sigma_{l}\rangle ^{1/2}\langle k_{j}\rangle}dk_{1}
d\tau_{1}dkd\tau \nonumber\\&&
\leq C\delta^{\frac{j}{2j+1}-2\epsilon}N^{-j}\|H\|_{L^{2}}\prod_{l=1}^{2}\|H_{l}\|_{L^{2}}.
\label{5.03}
\end{eqnarray}
We define $A=A_{1}\cup A_{2}\cup A_{3}$,
where
\begin{eqnarray*}
&&A=\left\{(k_{1},\tau_{1},k,\tau )\in \left(\dot{Z}\times \R\right)^{2}
:k=k_{1}+k_{2},\tau=\tau_{1}+\tau_{2},|k_{1}|\leq |k_{2}|,|k_{2}|\geq\frac{N}{2}\right\}\nonumber\\&&
A_{1}=\left\{(k_{1},\tau_{1},k,\tau )\in A:|k_{1}|\ll |k_{2}|,|k_{1}|\leq N\right\}\nonumber\\&&
A_{2}=\left\{(k_{1},\tau_{1},k,\tau )\in A:|k_{1}|\ll |k_{2}|,|k_{1}|> N\right\}\nonumber\\&&
A_{3}=\left\{(k_{1},\tau_{1},k,\tau )\in A: |k_{1}|\sim |k_{2}|\right\}.
\end{eqnarray*}
The integrals corrsponding to $A_{j}(j=1,2,3)$ will be denoted by $I_{1},I_{2},I_{3}.$
We consider cases
\begin{eqnarray*}
&&(a):\quad |\sigma|={\rm max}
\left\{|\sigma|,|\sigma_{1}|,|\sigma_{2}|\right\}
\geq C|k_{\rm min}||k_{\rm max}|^{2j},\nonumber\\
&&(b):\quad |\sigma_{1}|={\rm max}
\left\{|\sigma|,|\sigma_{1}|,|\sigma_{2}|\right\}
\geq C|k_{\rm min}||k_{\rm max}|^{2j},\nonumber\\
&&(c):\quad |\sigma_{2}|={\rm max}
\left\{|\sigma|,|\sigma_{1}|,|\sigma_{2}|\right\}
\geq C|k_{\rm min}||k_{\rm max}|^{2j}.\nonumber
\end{eqnarray*}
1. Estimate of $I_{1}.$
By using the mean value Theorem, we have that
\begin{eqnarray*}
m(k_{1}+k_{2})-m(k_{1})m(k_{2})=m^{\prime}(\theta k_{1}+k_{2})k_{1},
\end{eqnarray*}
thus in region $A_{1}$, we have that
$|\theta k_{1}+k_{2}|\sim |k_{2}|$ which yields that
\begin{eqnarray}
&&\left|\frac{m(k_{1}+k_{2})-m(k_{1})m(k_{2})}{m(k_{1})m(k_{2})}\right|=\frac{\left|m(k_{1}+k_{2})-m(k_{2})\right|}{m(k_{2})}\nonumber\\
&&\leq \frac{m^{\prime}(\theta k_{1}+k_{2})|k_{1}|}{m(k_{2})}\leq \frac{C|k_{1}|}{|k_{2}|}.\label{5.04}
\end{eqnarray}
When $(a)$ is valid, by using (\ref{5.04}), the Plancherel identity and H\"older inequality as well as  Lemma 2.11,
 we have that in this case the left hand side of  (\ref{5.03}) can be bounded by
\begin{eqnarray*}
&&\int_{\!\!\!\mbox{\scriptsize $
\begin{array}{l}
k=k_{1}+k_{2}\\
\tau=\tau_{1}+\tau_{2}
\end{array}
$}}\frac{|k_{1}||k|^{3}H(k,\tau)\prod_{l=1}^{2}H_{l}(k_{l},\tau_{l})}{|k_{2}|\langle\sigma\rangle ^{1/2}\langle k\rangle
\prod_{l=1}^{2}\langle\sigma_{l}\rangle ^{1/2}\langle k_{l}\rangle}dk_{1}d\tau_{1}
dkd\tau \nonumber\\&&\leq C\int_{\!\!\!\mbox{\scriptsize $
\begin{array}{l}
k=k_{1}+k_{2}\\
\tau=\tau_{1}+\tau_{2}
\end{array}
$}}\frac{|k|^{-j}|k_{1}|^{-1/2}H(k,\tau)\prod_{l=1}^{2}H_{l}(k_{l},\tau_{l})}{\langle k\rangle
\prod_{l=1}^{2}\langle\sigma_{l}\rangle ^{1/2}\langle k_{l}\rangle}dk_{1}d\tau_{1}
dkd\tau \nonumber
\\&&\leq CN^{-j}\|H\|_{L^{2}}\prod_{l=1}^{2}\left\|\mathscr{F}^{-1}
\left(\frac{H_{l}}{\langle\sigma_{l}\rangle^{1/2}}\right)\right\|_{L_{xt}^{4}}
\leq CN^{-j}\delta^{\frac{j}{2j+1}-2\epsilon}\|H\|_{L^{2}}\prod_{l=1}^{2}\|H_{l}\|_{L^{2}}.
\end{eqnarray*}
When $(b)$ is valid, by using (\ref{5.04}), the Plancherel identity and H\"older inequality as well as  Lemma 2.12,
 we have that in this case the left hand side of  (\ref{5.03}) can be bounded by
\begin{eqnarray*}
&&\int_{\!\!\!\mbox{\scriptsize $
\begin{array}{l}
k=k_{1}+k_{2}\\
\tau=\tau_{1}+\tau_{2}
\end{array}
$}}\frac{|k_{1}||k|^{3}H(k,\tau)\prod_{l=1}^{2}H_{l}(k_{l},\tau_{l})}{|k_{2}|\langle\sigma\rangle ^{1/2}\langle k\rangle
\prod_{l=1}^{2}\langle\sigma_{l}\rangle ^{1/2}\langle k_{l}\rangle}dk_{1}d\tau_{1}
dkd\tau \nonumber\\&&\leq C\int_{\!\!\!\mbox{\scriptsize $
\begin{array}{l}
k=k_{1}+k_{2}\\
\tau=\tau_{1}+\tau_{2}
\end{array}
$}}\frac{|k|^{-j}|k_{1}|^{-1/2}H(k,\tau)\prod_{l=1}^{2}H_{l}(k_{l},\tau_{l})}{\langle\sigma_{2}\rangle ^{1/2}
\langle\sigma\rangle ^{1/2}}dk_{1}d\tau_{1}dkd\tau \nonumber
\\&&\leq CN^{-j}\int_{\!\!\!\mbox{\scriptsize $
\begin{array}{l}
k=k_{1}+k_{2}\\
\tau=\tau_{1}+\tau_{2}
\end{array}
$}}\langle\sigma\rangle ^{-\frac{1}{2}+\epsilon}\frac{H(k,\tau)\prod_{l=1}^{2}H_{l}(k_{l},\tau_{l})}{\langle\sigma_{2}\rangle ^{1/2}
}dk_{1}d\tau_{1}dkd\tau \nonumber
\\
&&\leq  CN^{-j}
\left\|\int_{\!\!\!\mbox{\scriptsize $
\begin{array}{l}
k=k_{1}+k_{2}\\
\tau=\tau_{1}+\tau_{2}
\end{array}
$}}\langle\sigma\rangle ^{-\frac{1}{2}+\epsilon}\frac{\prod_{l=1}^{2}H_{l}(k_{l},\tau_{l})}{\langle\sigma_{2}\rangle ^{1/2}
}dk_{1}d\tau_{1}\right\|_{L^{2}}\|H\|_{L^{2}}\nonumber\\
&&\leq CN^{-j}\delta^{\frac{j}{2j+1}-2\epsilon}\|H\|_{L^{2}}\prod_{l=1}^{2}\|H_{l}\|_{L^{2}}.
\end{eqnarray*}
When $(c)$ is valid,  this case can be proved similarly to case $(b)$.

2. Estimate of $I_{2}.$
In this case, we have that
\begin{eqnarray*}
&&\frac{|m(k_{1}+k_{2})-m(k_{1})m(k_{2})|}{m(k_{1})m(k_{2})}\leq \frac{{\rm max}
\left\{m(k_{1}+k_{2}),m(k_{2})\right\}}{m(k_{1})m(k_{2})}\nonumber\\
&&\leq \frac{C}{m(k_{1})}\leq C \left(\frac{|k_{1}|}{N}\right)^{-s}.
\end{eqnarray*}
When $(a)$ is valid,  we have that  in this case the left hand side of  (\ref{5.03}) can be bounded by
\begin{eqnarray}
&&\int_{\!\!\!\mbox{\scriptsize $
\begin{array}{l}
k=k_{1}+k_{2}\\
\tau=\tau_{1}+\tau_{2}
\end{array}
$}}\frac{|k_{1}|^{-s}|k|^{3}N^{s}H(k,\tau)\prod_{l=1}^{2}H_{l}(k_{l},\tau_{l})}{\langle\sigma\rangle ^{1/2}\langle k\rangle
\prod_{l=1}^{2}\langle\sigma_{l}\rangle ^{1/2}\langle k_{l}\rangle}dk_{1}d\tau_{1}dkd\tau
\nonumber\\&&\leq C\int_{\!\!\!\mbox{\scriptsize $
\begin{array}{l}
k=k_{1}+k_{2}\\
\tau=\tau_{1}+\tau_{2}
\end{array}
$}}\frac{|k_{1}|^{-s-\frac{3}{2}}N^{s}|k|^{1-j}H(k,\tau)\prod_{l=1}^{2}H_{l}(k_{l},\tau_{l})}{
\prod_{l=1}^{2}\langle\sigma_{l}\rangle ^{1/2}}dk_{1}d\tau_{1}dkd\tau,\label{5.05}
\end{eqnarray}
if $-s-\frac{3}{2}\leq 0,$ by using the Plancherel identity and the H\"older inequality
as well as Lemma 2.11, we have that  (\ref{5.05})  can be bounded by
\begin{eqnarray*}
&&C\int_{\!\!\!\mbox{\scriptsize $
\begin{array}{l}
k=k_{1}+k_{2}\\
\tau=\tau_{1}+\tau_{2}
\end{array}
$}}\frac{N^{-s-\frac{3}{2}}N^{s}N^{1-j}H(k,\tau)\prod_{l=1}^{2}H_{l}(k_{l},\tau_{l})}{
\prod_{l=1}^{2}\langle\sigma_{l}\rangle ^{1/2}}dk_{1}d\tau_{1}dkd\tau\nonumber\\
&&\leq CN^{-j-\frac{1}{2}}\int_{\!\!\!\mbox{\scriptsize $
\begin{array}{l}
k=k_{1}+k_{2}\\
\tau=\tau_{1}+\tau_{2}
\end{array}
$}}\frac{H(k,\tau)\prod_{l=1}^{2}H_{l}(k_{l},\tau_{l})}{\prod_{l=1}^{2}\langle\sigma_{l}\rangle ^{1/2}}dk_{1}d\tau_{1}dkd\tau
\nonumber\\
&&\leq CN^{-j-\frac{1}{2}}\|H\|_{L^{2}}\prod_{l=1}^{2}\left\|\mathscr{F}^{-1}
\left(\frac{H_{l}}{\langle\sigma_{l}\rangle^{1/2}}\right)\right\|_{L_{xt}^{4}}\nonumber\\
&&\leq CN^{-j-\frac{1}{2}}\delta^{\frac{j}{2j+1}-2\epsilon}\|H\|_{L^{2}}\prod_{l=1}^{2}\|H_{l}\|_{L^{2}}.
\end{eqnarray*}
if $-s-\frac{3}{2}\geq 0$, since $s\geq \frac{2-j}{2}$, by using the Plancherel identity and the H\"older inequality
as well as Lemma 2.11, we have that  (\ref{5.05})  can be bounded by
\begin{eqnarray*}
&&C\int_{\!\!\!\mbox{\scriptsize $
\begin{array}{l}
k=k_{1}+k_{2}\\
\tau=\tau_{1}+\tau_{2}
\end{array}
$}}\frac{|k|^{-s-\frac{3}{2}}N^{s}|k|^{1-j}H(k,\tau)\prod_{l=1}^{2}H_{l}(k_{l},\tau_{l})}{
\prod_{l=1}^{2}\langle\sigma_{l}\rangle ^{1/2}}dk_{1}d\tau_{1}dkd\tau\nonumber\\
&&\leq C\int_{\!\!\!\mbox{\scriptsize $
\begin{array}{l}
k=k_{1}+k_{2}\\
\tau=\tau_{1}+\tau_{2}
\end{array}
$}}\frac{|k|^{-s-\frac{1}{2}-j}N^{s}H(k,\tau)\prod_{l=1}^{2}
H_{l}(k_{l},\tau_{l})}{\prod_{l=1}^{2}\langle\sigma_{l}\rangle ^{1/2}}dk_{1}d\tau_{1}dkd\tau
\nonumber\\
&&\leq CN^{-j-\frac{1}{2}}\int_{\!\!\!\mbox{\scriptsize $
\begin{array}{l}
k=k_{1}+k_{2}\\
\tau=\tau_{1}+\tau_{2}
\end{array}
$}}\frac{H(k,\tau)\prod_{l=1}^{2}H_{l}(k_{l},\tau_{l})}{\prod_{l=1}^{2}
\langle\sigma_{l}\rangle ^{1/2}}dk_{1}d\tau_{1}dkd\tau\nonumber\\
&&\leq CN^{-j-\frac{1}{2}}\|H\|_{L^{2}}\prod_{l=1}^{2}\left\|\mathscr{F}^{-1}
\left(\frac{H_{l}}{\langle\sigma_{l}\rangle^{1/2}}\right)\right\|_{L_{xt}^{4}}
\leq CN^{-j-\frac{1}{2}}\delta^{\frac{j}{2j+1}-2\epsilon}\|H\|_{L^{2}}\prod_{l=1}^{2}\|H_{l}\|_{L^{2}}.
\end{eqnarray*}
When $(b)$ is valid,  by using (\ref{5.04})  and the Plancherel identity and the H\"older inequality
as well as Lemma 2.11,  we have that in this case the left hand side of (\ref{5.03}) can be bounded by
\begin{eqnarray}
&&\int_{\!\!\!\mbox{\scriptsize $
\begin{array}{l}
k=k_{1}+k_{2}\\
\tau=\tau_{1}+\tau_{2}
\end{array}
$}}\frac{|k_{1}|^{-s}|k|^{3}N^{s}H(k,\tau)\prod_{l=1}^{2}H_{l}(k_{l},\tau_{l})}{\langle\sigma\rangle ^{1/2}\langle k\rangle
\prod_{l=1}^{2}\langle\sigma_{l}\rangle ^{1/2}\langle k_{l}\rangle}dk_{1}d\tau_{1}dkd\tau
\nonumber\\&&\leq C\int_{\!\!\!\mbox{\scriptsize $
\begin{array}{l}
k=k_{1}+k_{2}\\
\tau=\tau_{1}+\tau_{2}
\end{array}
$}}\frac{|k_{1}|^{-s-\frac{3}{2}}N^{s}|k|^{1-j}H\prod_{l=1}^{2}H_{l}}{\langle\sigma\rangle ^{1/2}
\langle\sigma_{2}\rangle ^{1/2}}dk_{1}d\tau_{1}dkd\tau,\label{5.06}
\end{eqnarray}
if $-s-\frac{3}{2}\leq 0,$  by using the Plancherel identity and the H\"older inequality
as well as Lemma 2.12,   we have that  (\ref{5.06})  can be bounded by
\begin{eqnarray*}
&&C\int_{\!\!\!\mbox{\scriptsize $
\begin{array}{l}
k=k_{1}+k_{2}\\
\tau=\tau_{1}+\tau_{2}
\end{array}
$}}\frac{N^{-s-\frac{3}{2}}N^{s}N^{1-j}H(k,\tau)\prod_{l=1}^{2}H_{l}(k_{l},\tau_{l})}{\langle\sigma\rangle ^{1/2}
\langle\sigma_{2}\rangle ^{1/2}}dk_{1}d\tau_{1}dkd\tau\nonumber\\
&&\leq CN^{-j-\frac{1}{2}}\int_{\!\!\!\mbox{\scriptsize $
\begin{array}{l}
k=k_{1}+k_{2}\\
\tau=\tau_{1}+\tau_{2}
\end{array}
$}}\frac{H(k,\tau)\prod_{l=1}^{2}H_{l}(k_{l},\tau_{l})}{\langle\sigma\rangle ^{1/2}
\langle\sigma_{2}\rangle ^{1/2}}dk_{1}d\tau_{1}dkd\tau\nonumber\\
&&\leq CN^{-j-\frac{1}{2}}\left\|\langle \sigma\rangle ^{-\frac{1}{2}+\epsilon}\int_{\!\!\!\mbox{\scriptsize $
\begin{array}{l}
k=k_{1}+k_{2}\\
\tau=\tau_{1}+\tau_{2}
\end{array}
$}}\frac{\prod_{l=1}^{2}H_{l}(k_{l},\tau_{l})}{
\langle\sigma_{2}\rangle ^{1/2}}dk_{1}d\tau_{1} \right\|_{L^{2}}\|H\|_{L^{2}}\nonumber\\
&&\leq CN^{-j-\frac{1}{2}}\delta^{\frac{j}{2j+1}-2\epsilon}\|H\|_{L^{2}}\prod_{l=1}^{2}\|H_{l}\|_{L^{2}}.
\end{eqnarray*}
if $-s-\frac{3}{2}\geq 0$, since $s\geq \frac{2-j}{2}$,  (\ref{5.06})  can be bounded by
\begin{eqnarray*}
&&C\int_{\!\!\!\mbox{\scriptsize $
\begin{array}{l}
k=k_{1}+k_{2}\\
\tau=\tau_{1}+\tau_{2}
\end{array}
$}}\frac{|k|^{-s-\frac{3}{2}}N^{s}|k|^{1-j}H(k,\tau)\prod_{l=1}^{2}H_{l}(k_{l},\tau_{l})}{\langle\sigma\rangle^{1/2}
\prod_{l=1}^{2}\langle\sigma_{l}\rangle ^{1/2}}dk_{1}d\tau_{1}dkd\tau\nonumber\\
&&\leq C\int_{\!\!\!\mbox{\scriptsize $
\begin{array}{l}
k=k_{1}+k_{2}\\
\tau=\tau_{1}+\tau_{2}
\end{array}
$}}\frac{|k|^{-s-\frac{1}{2}-j}N^{s}H(k,\tau)\prod_{l=1}^{2}H_{l}(k_{l},\tau_{l})}
{\langle\sigma\rangle^{1/2}\langle\sigma_{1}\rangle ^{1/2}}dk_{1}d\tau_{1}dkd\tau
\nonumber\\
&&\leq CN^{-j-\frac{1}{2}}\int_{\!\!\!\mbox{\scriptsize $
\begin{array}{l}
k=k_{1}+k_{2}\\
\tau=\tau_{1}+\tau_{2}
\end{array}
$}}\frac{H(k,\tau)\prod_{l=1}^{2}H_{l}(k_{l},\tau_{l})}{\langle\sigma\rangle ^{1/2}\langle\sigma_{2}\rangle ^{1/2}}dk_{1}d\tau_{1}dkd\tau
\nonumber\\
&&\leq CN^{-j-\frac{1}{2}}\left\|\langle \sigma\rangle ^{-\frac{1}{2}+\epsilon}\int_{\!\!\!\mbox{\scriptsize $
\begin{array}{l}
k=k_{1}+k_{2}\\
\tau=\tau_{1}+\tau_{2}
\end{array}
$}}\frac{\prod_{l=1}^{2}H_{l}}{
\langle\sigma_{2}\rangle ^{1/2}}dk_{1}d\tau_{1} \right\|_{L^{2}}\|H\|_{L^{2}}\nonumber\\
&&\leq CN^{-j-\frac{1}{2}}\delta^{\frac{j}{2j+1}-2\epsilon}\|H\|_{L^{2}}\prod_{l=1}^{2}\|H_{l}\|_{L^{2}}.
\end{eqnarray*}
When $(c)$ is valid,  this case can be proved similarly to case $(b)$.

3. Estimate of $I_{3}.$
In this case, we have that
\begin{eqnarray}
&&\frac{|m(k_{1}+k_{2})-m(k_{1})m(k_{2})|}{m(k_{1})m(k_{2})}\leq
C\prod_{l=1}^{2} \left(\frac{|k_{l}|}{N}\right)^{-s}.\label{5.07}
\end{eqnarray}
When $(a)$ is valid, by using (\ref{5.07}) and the Plancherel identity and the H\"older inequality
as well as Lemma 2.11,  since $\frac{2-j}{2}\leq s\leq 1,$  we have that in this case the left hand side of (\ref{5.03}) can be bounded by
\begin{eqnarray*}
&&\int_{\!\!\!\mbox{\scriptsize $
\begin{array}{l}
k=k_{1}+k_{2}\\
\tau=\tau_{1}+\tau_{2}
\end{array}
$}}\frac{|k|^{3}|k_{1}|^{-2s}N^{2s}H(k,\tau)\prod_{l=1}H_{l}(k_{l},\tau_{l})}{\langle\sigma\rangle ^{1/2}\langle k\rangle
\prod_{l=1}^{2}\langle\sigma_{l}\rangle ^{1/2}\langle k_{l}\rangle}dk_{1}d\tau_{1}dkd\tau
\nonumber\\&&\leq C\int_{\!\!\!\mbox{\scriptsize $
\begin{array}{l}
k=k_{1}+k_{2}\\
\tau=\tau_{1}+\tau_{2}
\end{array}
$}}\frac{|k_{1}|^{-2s-2-j}N^{2s}|k|^{5/2}H(k,\tau)\prod_{l=1}^{2}H_{l}(k_{l},\tau_{l})}{\langle k\rangle
\prod_{l=1}^{2}\langle\sigma_{l}\rangle ^{1/2}}
dk_{1}d\tau_{1}dkd\tau \nonumber\\&&
\leq C\int_{\!\!\!\mbox{\scriptsize $
\begin{array}{l}
k=k_{1}+k_{2}\\
\tau=\tau_{1}+\tau_{2}
\end{array}
$}}\frac{|k_{1}|^{-2s-\frac{1}{2}-j}N^{2s}H(k,\tau)\prod_{l=1}^{2}H_{l}(k_{l},\tau_{l})}{
\prod_{l=1}^{2}\langle\sigma_{l}\rangle ^{1/2}}
dk_{1}d\tau_{1}dkd\tau \nonumber\\&&
\leq CN^{-j-\frac{1}{2}}\int_{\!\!\!\mbox{\scriptsize $
\begin{array}{l}
k=k_{1}+k_{2}\\
\tau=\tau_{1}+\tau_{2}
\end{array}
$}}\frac{H(k,\tau)\prod_{l=1}^{2}H_{l}(k_{l},\tau_{l})}{
\prod_{l=1}^{2}\langle\sigma_{l}\rangle ^{1/2}}
dk_{1}d\tau_{1}dkd\tau \nonumber\\&&
\leq CN^{-j-\frac{1}{2}}
\|H\|_{L^{2}}\prod_{l=1}^{2}\left\|\mathscr{F}^{-1}\left(\frac{H_{l}}{\langle\sigma_{l}
\rangle^{1/2}}\right)\right\|_{L_{xt}^{4}}\nonumber\\
&&\leq CN^{-j-\frac{1}{2}}\delta^{\frac{j}{2j+1}-2\epsilon}\|H\|_{L^{2}}\prod_{l=1}^{2}\|H_{l}\|_{L^{2}}.
\end{eqnarray*}
When $(b)$ is valid, by using (\ref{5.07}) and the Plancherel identity and the H\"older inequality
as well as Lemma 2.12,  since $\frac{2-j}{2}\leq s\leq 1,$  we have that in this case the left hand side of (\ref{5.03}) can be bounded by
\begin{eqnarray*}
&&\int_{\!\!\!\mbox{\scriptsize $
\begin{array}{l}
k=k_{1}+k_{2}\\
\tau=\tau_{1}+\tau_{2}
\end{array}
$}}\frac{|k|^{3}|k_{1}|^{-2s}N^{2s}H(k,\tau)\prod_{l=1}H_{l}(k_{l},\tau_{l})}{\langle\sigma\rangle ^{1/2}\langle k\rangle
\langle\sigma_{2}\rangle ^{1/2}\prod_{l=1}^{2}\langle k_{l}\rangle}dk_{1}d\tau_{1}dkd\tau
\nonumber\\&&\leq C\int_{\!\!\!\mbox{\scriptsize $
\begin{array}{l}
k=k_{1}+k_{2}\\
\tau=\tau_{1}+\tau_{2}
\end{array}
$}}\frac{|k_{1}|^{-2s-2-j}N^{2s}|k|^{5/2}H(k,\tau)\prod_{l=1}^{2}H_{l}(k_{l}\tau_{l})}{\langle k\rangle\langle\sigma\rangle ^{1/2}
\langle\sigma_{2}\rangle ^{1/2}}
dk_{1}d\tau_{1}dkd\tau \nonumber\\&&
\leq C\int_{\!\!\!\mbox{\scriptsize $
\begin{array}{l}
k=k_{1}+k_{2}\\
\tau=\tau_{1}+\tau_{2}
\end{array}
$}}\frac{|k_{1}|^{-2s-\frac{1}{2}-j}N^{2s}H(k,\tau)\prod_{l=1}^{2}H_{l}(k_{l},\tau_{l})}{\langle\sigma\rangle ^{1/2}
\langle\sigma_{2}\rangle ^{1/2}}
dk_{1}d\tau_{1}dkd\tau \nonumber\\&&
\leq CN^{-j-\frac{1}{2}}\int_{\!\!\!\mbox{\scriptsize $
\begin{array}{l}
k=k_{1}+k_{2}\\
\tau=\tau_{1}+\tau_{2}
\end{array}
$}}\frac{H(k,\tau)\prod_{l=1}^{2}H_{l}(k_{l},\tau_{l})}{\langle\sigma\rangle ^{1/2}
\langle\sigma_{2}\rangle ^{1/2}}
dk_{1}d\tau_{1}dkd\tau \nonumber\\&&
\leq CN^{-j-\frac{1}{2}}\left\|\langle \sigma\rangle ^{-\frac{1}{2}+\epsilon}\int_{\!\!\!\mbox{\scriptsize $
\begin{array}{l}
k=k_{1}+k_{2}\\
\tau=\tau_{1}+\tau_{2}
\end{array}
$}}\frac{\prod_{l=1}^{2}H_{l}(k_{l},\tau_{l})}{
\langle\sigma_{2}\rangle ^{1/2}}dk_{1}d\tau_{1}\right\|_{L^{2}}\|H\|_{L^{2}}\nonumber\\
&&\leq CN^{-j-\frac{1}{2}}\delta^{\frac{j}{2j+1}-2\epsilon}\|H\|_{L^{2}}\prod_{l=1}^{2}\|H_{l}\|_{L^{2}}.
\end{eqnarray*}
When $(c)$ is valid,  this case can be proved similarly to case $(b)$.

We have completed the proof of Lemma 5.1.

\begin{Lemma}\label{Lemma5.2}
Let $\frac{j-2}{2}\leq s<1$ and $u$ be the solution to the Cauchy problem for (\ref{1.01}) on $[0,\delta]$. Then
\begin{eqnarray}
\left|\int_{0}^{\delta}\int_{\mathbf{T}}(\partial_{x}(Iu)\left[I(u_{x}^{2})-(\partial_{x}Iu)^{2}\right]dxdt\right|
\leq C\delta^{\frac{j}{2j+1}-2\epsilon}N^{-j}\|Iu\|_{X_{1,\frac{1}{2}}^{\delta}}^{3}\label{5.08}.
\end{eqnarray}
\end{Lemma}

Lemma 5.2 can be proved similarly to Lemma 5.1.

\begin{Lemma}\label{Lemma5.3}
Let $-\frac{j}{2}\leq s<1$ and $u$ be the solution to the Cauchy problem for (\ref{1.01}) on $[0,\delta]$. Then
\begin{eqnarray}
\left|\int_{0}^{\delta}\int_{\mathbf{T}}\partial_{x}(Iu)\left[Iu^{2}-(Iu)^{2}\right]dxdt\right|
\leq C\delta^{\frac{j}{2j+1}-2\epsilon}N^{-j-2}\|Iu\|_{X_{1,\frac{1}{2}}^{\delta}}^{3}\label{5.09}.
\end{eqnarray}
\end{Lemma}

Lemma 5.3 can be proved similarly to Lemma 5.1.
\begin{Lemma}\label{Lemma5.4}
Let $\frac{2-j}{2}\leq s<1$ and $u$ be the solution to the Cauchy problem for (\ref{1.01})
on $[0,\delta]$. Then
\begin{eqnarray}
      \left|\|Iu(\delta)\|_{H^{1}}^{2}-\|Iu(0)\|_{H^{1}}^{2}\right|\leq
      C\delta^{\frac{j}{2j+1}-2\epsilon}N^{-j}\|Iu\|_{X_{1,\frac{1}{2}}^{\delta}}^{3}\label{5.010}
\end{eqnarray}
\end{Lemma}
{\bf Proof.} By using a proof similar to  (4.3) of \cite{LY}, we have that
\begin{eqnarray}
&&\|Iu(\delta)\|_{H^{1}}^{2}-\|Iu(0)\|_{H^{1}}^{2}=\int_{0}^{\delta}
\int_{\mathbf{T}}(1-\partial_{x}^{2})\partial_{x}(Iu)\left[Iu^{2}-(Iu)^{2}\right]dxdt
\nonumber\\&&\quad+2\int_{0}^{\delta}
\int_{\mathbf{T}}(\partial_{x}(Iu)\left[Iu^{2}-(Iu)^{2}\right]dxdt
\nonumber\\&&\quad+\int_{0}^{\delta}\int_{\mathbf{T}}(\partial_{x}(Iu)\left[I(u_{x}^{2})-(\partial_{x}Iu)^{2}\right]
dxdt\label{5.011}
\end{eqnarray}
{\bf Proof.} To prove (\ref{5.011}), it suffices to prove that
\begin{eqnarray}
&&\left|\|Iu(\delta)\|_{H^{1}}^{2}-\|Iu(0)\|_{H^{1}}^{2}\right|\leq\left|\int_{0}^{\delta}
\int_{\mathbf{T}}(1-\partial_{x}^{2})\partial_{x}(Iu)\left[Iu^{2}-(Iu)^{2}\right]dxdt\right|
\nonumber\\&&\quad+2\int_{0}^{\delta}
\left|\int_{\mathbf{T}}(\partial_{x}(Iu)\left[Iu^{2}-(Iu)^{2}\right]dxdt\right|\nonumber\\
&&\quad+\left|\int_{\mathbf{T}}(\partial_{x}(Iu)\left[I(u_{x}^{2})
-(\partial_{x}Iu)^{2}\right]dxdt\right|\leq C\delta^{\frac{j}{2j+1}-2\epsilon}N^{-j}
\|Iu\|_{X_{1,\frac{1}{2}}^{\delta}}^{3}\label{5.012}.
\end{eqnarray}
 (\ref{5.012}) can be obtained from Lemmas 5.1-5.3.

We have completed the proof of Lemma 5.4.

\bigskip
\bigskip

\noindent {\large\bf 6. Proof of Theorem  1.2}

\setcounter{equation}{0}

 \setcounter{Theorem}{0}

\setcounter{Lemma}{0}

\setcounter{section}{6}

We give Theorem 5.1 which is a variant of Theorem 1.1 before giving the proof of Theorem 1.2.

We consider the Cauchy problem for
\begin{eqnarray}
&&(Iu)_{t}+\partial_{x}^{2j+1}(Iu)+\frac{1}{2}\partial_{x}I(u^{2})+\partial_{x}
(1-\partial_{x}^{2})^{-1}I\left[u^{2}+\frac{1}{2}u_{x}^{2}\right]=0,\label{6.01}\\
&&Iu(x,0)=Iu_{0}(x).\label{6.02}
\end{eqnarray}
\begin{Theorem}\label{Thm5.1}
Let $s\geq -\frac{j-2}{2}$  and $u_{0}$ be $2\pi$-periodic
function and  zero $x$-mean and $Iu_{0}\in H^{1}(\mathbf{T})$.
Then the Cauchy problems (\ref{6.01})(\ref{6.02})
are locally well-posed.
\end{Theorem}
{\bf Proof.} Let $Iu=v$, we define
\begin{eqnarray*}
&&G(v)=\eta(t)S(t)v(0)\nonumber\\&&+\eta(t)\int_{0}^{t}\left[\frac{1}{2}\partial_{x}
I(\eta\left(\frac{t}{\delta}\right)u)^{2}
+\partial_{x}(1-\partial_{x}^{2})^{-1}I\left[(\eta\left(\frac{t}{\delta}\right)u)^{2}
+\frac{1}{2}(\eta\left(\frac{t}{\delta}\right)u_{x})^{2}\right]\right]dt^{'}.
\end{eqnarray*}
and
\begin{eqnarray}
B=\left\{u\in Y_{1}^{\delta}: \quad \|Iu\|_{ Y_{1}^{\delta}}\leq 2C\|Iu_{0}\|_{H^{1}
(\mathbf{T})}\right\},\label{6.03}
\end{eqnarray}
and
\begin{eqnarray*}
&&E=\frac{1}{2}\partial_{x}I(\eta\left(\frac{t}{\delta}\right)u)^{2}
+\partial_{x}(1-\partial_{x}^{2})^{-1}I\left[(\eta\left(\frac{t}{\delta}\right)u)^{2}
+\frac{1}{2}(\eta\left(\frac{t}{\delta}\right)u_{x})^{2}\right]\nonumber\\&&\qquad-
\frac{1}{2}\partial_{x}(\eta\left(\frac{t}{\delta}\right)Iv)^{2}
-\partial_{x}(1-\partial_{x}^{2})^{-1}\left[(\eta\left(\frac{t}{\delta}\right)Iv)^{2}
+\frac{1}{2}(\eta\left(\frac{t}{\delta}\right)Iv_{x})^{2}\right].
\end{eqnarray*}
Thus, we have that
\begin{eqnarray*}
&&G(v)=\eta(t)S(t)Iu_{0}\nonumber\\&&+\eta(t)\int_{0}^{t}\left[E+\partial_{x}
(1-\partial_{x}^{2})^{-1}\left[(\eta\left(\frac{t}{\delta}\right)Iv)^{2}+\frac{1}{2}
(\eta\left(\frac{t}{\delta}\right)Iv_{x})^{2}\right]\right]dt^{'}.
\end{eqnarray*}
By using Lemmas 3.7-3.9, 5.1-5.3,  we have that
\begin{eqnarray*}
&&\|G(v)\|_{Y_{1}^{\delta}}\nonumber\\&&\leq \left\|\eta(t)S(t)Iu_{0}\right\|_{Y_{1}^{\delta}}
+\left\|\eta(t)\int_{0}^{t}\left[\partial_{x}
(1-\partial_{x}^{2})^{-1}\left[(\eta\left(\frac{t}{\delta}\right)Iv)^{2}+\frac{1}{2}
(\eta\left(\frac{t}{\delta}\right)Iv_{x})^{2}\right]\right]dt^{'}\right\|_{Y_{1}^{\delta}}\nonumber\\&&
\quad +\left\|\eta(t)\int_{0}^{t}Edt^{'}\right\|_{Y_{1}^{\delta}}\nonumber\\&&\leq
C\|Iu_{0}\|_{H^{1}}+C\delta^{\frac{j}{2j+1}-3\epsilon}\left\|v\right\|_{Y_{1}^{\delta}}^{2}\leq
2C\|Iu_{0}\|_{H^{1}}.
\end{eqnarray*}
Thus, $G$ maps $B$ into $B$.
By using Lemmas 3.7-3.9, 5.1-5.3,  we have that
\begin{eqnarray*}
\left\|G(u)-G(v)\right\|_{Y_{1}^{\delta}}\leq \frac{1}{2}\left\|u-v\right\|_{Y_{1}^{\delta}}.
\end{eqnarray*}
$G$  is a contraction mapping.

We have completed the proof of Theorem 5.1.

Now we are in a position to prove Theorem 1.2.
For $u_{0}\in H^{s}(\mathbf{T})$, from Theorem 5.1, we have that $u$ exists on $[0,\delta]$
and
\begin{eqnarray}
\delta \sim \|Iu_{0}\|_{H^{1}}^{-\frac{2j+1}{j-3(2j+1)\epsilon}}.\label{6.04}
\end{eqnarray}
From Theorem 5.1, we have that
\begin{eqnarray}
\|Iu\|_{Y_{1}^{\delta}}\leq 2C\|Iu_{0}\|_{H^{1}}.\label{6.05}
\end{eqnarray}
Combining (\ref{6.05}) with Lemma 5.4, we have that
\begin{eqnarray}
\|Iu(\delta)\|_{H^{1}}^{2}\leq \|Iu_{0}\|_{H^{1}}^{2}+CN^{-j}\delta^{\frac{j}{2j+1}-3\epsilon}\|Iu_{0}\|_{H^{1}}^{3}.\label{6.06}
\end{eqnarray}
If
\begin{eqnarray}
CN^{-j}\delta^{\frac{j}{2j+1}-3\epsilon}\|Iu_{0}\|_{H^{1}}^{3}\leq 3\|Iu_{0}\|_{H^{1}}^{2},\label{6.07}
\end{eqnarray}
then, we have that
\begin{eqnarray}
\|Iu(\delta)\|_{H^{1}}\leq 2\|Iu_{0}\|_{H^{1}},\label{6.08}
\end{eqnarray}
thus, we can consider $u(\delta)$ as the initial data, repeat the above process and extend the local solution on $[0,\delta]$ to
the local solution on $[\delta,2\delta]$. To extend the local solution to the global on time interval $[0,T]$, we need to extend $[T\delta^{-1}]$ times,
from (\ref{6.07}), it suffices to prove that
\begin{eqnarray}
CN^{-j}\delta^{\frac{j}{2j+1}-3\epsilon}\|Iu_{0}\|_{H^{1}}^{3}T\delta^{-1}\leq 3\|Iu_{0}\|_{H^{1}}^{2},\label{6.09}
\end{eqnarray}
It is easily checked that
\begin{eqnarray}
\|u\|_{H^{s}}\leq \|Iu_{0}\|_{H^{1}}\leq CN^{1-s}\|u\|_{H^{s}}.\label{6.010}
\end{eqnarray}
Combining (\ref{6.04}), (\ref{6.010}) with (\ref{6.09}), we have that
\begin{eqnarray}
CTN^{\left[\frac{(2j+1)(1-s)}{j-3(2j+1)\epsilon}\right]
(1-s)-j}\|u_{0}\|_{H^{s}}^{\frac{2j+1}{j-3(2j+1)\epsilon}}\leq 1.\label{6.011}
\end{eqnarray}
Let $f(j)=\frac{(2j+1)}{j-3(2j+1)\epsilon}$.
To obtain (\ref{6.011}), it suffices to choose $s>\frac{2j+1-j^{2}}{2j+1}$ and
\begin{eqnarray}
N=\left(CT\|u_{0}\|_{H^{s}}^{f(j)}\right)^{\frac{1}{j-f(j)(1-s)}}.\label{6.012}
\end{eqnarray}
From the above iteration process, we have that
\begin{eqnarray}
&&\sup\limits_{t\in [0,T]}\|u(\cdot,t)\|_{H^{s}}\leq 2\|Iu_{0}\|_{H^{1}}\nonumber\\&&
\leq CN^{1-s}\|u_{0}\|_{H^{s}}\leq C\left(CT\|u_{0}\|_{H^{s}}^{f(j)}\right)^{\frac{1-s}{j-f(j)(1-s)}}\|u_{0}\|_{H^{s}}\nonumber\\
&&\leq CT^{\frac{1-s}{j-f(j)(1-s)}}\|u_{0}\|_{H^{s}}^{\frac{j}{j-f(j)(1-s)}}.
\end{eqnarray}

We have completed the proof of Theorem 1.2.

\leftline{\large \bf Acknowledgments}

\bigskip

\noindent

 This work is supported by the Natural Science Foundation of China
 under grant numbers 11171116 and 11401180. The second author is also
 supported in part by the Fundamental Research Funds for the
 Central Universities of China under the grant number 2012ZZ0072.
 The third author is  supported by the
 NSF of China (No.11371367) and Fundamental
 research program of NUDT(JC12-02-03).


\bigskip

  \bigskip


\leftline{\large\bf  References}



\begin{thebibliography}{99}






\bibitem{Bou} J. Bourgain, On the Cauchy problem for the Kadomtsev-Petviashvili equation, {\it Geom. Funct. Anal.}
3(1993), 115-159.

 \bibitem{B}
     J. Bourgain,
     Fourier transform restriction phenomena for certain lattice subsets and
     applications to nonlinear evolution equations,  part I: Schr\"odinger equations,  {\it Geom. Funct. Anal.} 3(1993),  107-156.
\bibitem{Bourgain-GAFA93}
        J. Bourgain,
        Fourier transform restriction phenomena for certain
        lattice subsets and applications to nonlinear evolution equations,
        part II: The KdV equation,
        {\it Geom.   Funct. Anal.,} 3(1993), 209-262.

\bibitem{Bou}
 J. Bourgain, Periodic Korteweg  de vries equation with measures as initial data, {\it  Sel. Math.} 3(1997), 115-159.

\bibitem{By} P.J. Byers, The initial value problem for a KdV-type equation and related bilinear estimate, Dissertation, University
of Notre Dame, 2003.


\bibitem{CH} R. Camassa, D. Holm, An integrable shallow water equation with peaked solutions, {\it Phys. Rev. Lett.,}  71(1993),
1661¨C1664.







\bibitem{CKSTT}  J. Colliander, M. Keel, G. Staffilani,  H. Takaoka,  T. Tao,
Sharp global well-posedness for KdV and modified KdV on $\R$  and $\mathbf{T}$,   {\it J. Amer. Math. Soc.} 16(2003),  705-749.



\bibitem{FF}A. Fokas, B. Fuchssteiner, Symplectic structures, their $B\ddot{a}$klund transformations and hereditary symmetries, {\it Phys.
Rev. Lett.,}  71(1981),  47-66.







\bibitem{G}
J. Gorsky, On the Cauchy problem for a KdV-type equation on the circle, {Notre  Dame: University of  Notre  Dame, 2004}.
\bibitem{G}
A. Gr$\ddot{u}$nrock, New applications of the Fourier restriction norm method to well-posedness problemds for nonlinear evolution equations,
{\it Wuppertal: University of Wuppertal, 2002}


\bibitem{HM1998} A. A. Himonas, G. Misiolek, The  Cauchy   problem for  a shallow water type equation, {\it  Comm. Partial   Diff. Eqns.}
23(1998), 123-139.




\bibitem{HM2000} A. A. Himonas, G. Misiolek,  Well-posedness of the Cauchy problem for a shallow water equation on the circle,
{\it J. Diff. Eqns.} 161(2000), 479-495.



\bibitem{HM1} A.A. Himonas, G. Misiolek, The initial value problem for a fifth order shallow water, in: Analysis, Geometry,
Number Theory: The Mathematics of Leon Ehrenpreis, in: Contemp. Math., vol. 251, Amer. Math. Soc., Providence,
RI, 2000, pp. 309¨C320







\bibitem{H} H. Hirayama, LocaL well-posedness for the  periodic higher order KdV type equations,
 {\it Nonlinear Differential euqtaions and applications,}
19(2012), 677-693.





\bibitem{IK} A. D. Ionescu, C. E. Kenig, Global well-posedness of the Benjamin-Ono equation
in low-regularity spaces, {\it J. Amer. Math. Soc.} 20(2007), 753-798.

\bibitem{IKT} A. D. Ionescu, C. E. Kenig, D. Tataru,  Global well-posedness of the KP-I
initial-value problem in the energy space, {\it Invent. Math.} 173(2008), 265-304.

\bibitem{KT2005} T. Kappeler and P. Topalov, Global well-posedness of mKdV in
$L^2(T,R),${\it  Comm. Partial Differential Equations} 30(2005),
435-449.
\bibitem{KT2006} T. Kappeler and P. Topalov, Global wellposedness of KdV in
$H^{-1}(T,R)$, {\it Duke Math. J.}  135(2006),  327-360.


\bibitem{LJ}X. Liu, Y. Jin, The Cauchy problem of a shallow water equation, {\it Acta Math. Sin. (Engl. Ser.)} 30(2004), 1-16.





\bibitem{O}E. A. Olson,
Well posedness for a higher order modified Camassa-Holm equation,
{\it J.  Diff.  Eqns.}  246(2009),  4154-4172.

\bibitem{KPV}
C. E. Kenig,  G. Ponce, L. Vega, A bilinear estimate with applications to the KdV equation, {\it J. Amer. Math. Soc.} 9(1996), 573-603.

\bibitem{KPV0}
C. E. Kenig,  G. Ponce, L. Vega, On the ill-posedness of some canonical dispersive equations, {\it Duke Math. J.} 106(2001), 617-633.

\bibitem{Kis} N. Kishimoto, Well-posedness of the Cauchy problem for the Korteweg-de
Vries
equation at the critical regularity, {\it Diff. Int. Eqns.}  22(2009), 447-464.

\bibitem{LYY}Y. S. Li,  W. Yan,   X. Y.  Yang,
Well-posedness of a higher order modified Camassa-Holm equation in spaces of low regularity, {\it J. Evol. Eqns.} 10(2010),  465-486.

\bibitem{M} L. Molinet, A note on ill-posedness for the KdV equation, Differential
Integral Equations 24 (2011), 759¨C765.




\bibitem{LY} Y. S. Li, X. Y. Yang, Global well-posedness for a fifth-order shallow water equation on the circle, {\it Acta Mathematica Scientia}, 31(2011), 1303-1317.
\bibitem{Molinet} L. Molinet,  Sharp ill-posedness results for the KdV and mKdV
equations on the torus, {Advances in Mathematics,}  230(2012), 1895-1930.

\bibitem{ST} J. C. Saut, N. Tzvetov, On the periodic  KP-I type equations,
 {\it Comm. Math. Phys.} 221(2001), 451-476.

\bibitem{T}  T. Tao, Multilinear weighted convolution of $L^2$ functions, and
applications
to non-linear dispersive equations, {\it Amer. J. Math., } 123(2001), 839-908.


\bibitem{WC}  H. Wang, S. B. Cui,  Global well-posedness of the Cauchy problem
of the fifth-order shallow water equation, {\it J. Diff. Eqns.,} 230(2006), 600-613.








\bibitem{YLL} W. Yan, Y.S. Li, S.M. Li,  Sharp well-posedness and ill-posedness of a higher-order
modified Camassa-Holm equation, {\it  Diff. Int. Eqns.,} 25(2012), 1053-1074.
\bibitem{YL}X. Y. Yang,  Y. S.  Li, Global well-posedness for a fifth-order shallow water equation in Sobolev spaces, {\it J. Diff. Eqns.} 248(2010), 1458-1472.

\end{thebibliography}
\end{document}